\documentclass[11pt]{amsart}

\usepackage{amsmath,amsthm, amscd, amssymb, amsfonts}
\usepackage{epsfig}

\newcommand{\Hc}{{\mathcal H}}

\newcommand{\lag}{\mbox{\rm Lag\,}}
\newcommand{\otb}{{\overline{\otimes}}}
\newcommand{\otk}{{\otimes_{\ku}}}

\newcommand{\Mo}{{\mathcal M}}

\newcommand{\No}{{\mathcal N}}

\newcommand{\Kc}{{\mathcal K}}
\newcommand{\nic}{{\mathfrak B}}

\newcommand{\oc}{{\mathcal O}}
\newcommand{\ocb}{\overline{{\mathcal O}}}
\newcommand{\Ss}{{\mathcal S}}
\newcommand{\ot}{{\otimes}}

\newcommand{\kc}{{\mathcal K}}
\newcommand{\Ac}{{\mathcal A}}
\newcommand{\ere}{{\mathcal R}}
\newcommand{\ele}{{\mathcal L}}

\newcommand{\phib}{{\overline \phi}}
\newcommand{\ca}{{\mathcal C}}

\newcommand{\lagb}{{\overline \lag}}

\newcommand{\Do}{{\mathcal D}}

\newcommand{\Fc}{{\mathcal F}}
\newcommand{\Gc}{{\mathcal G}}

\newcommand{\op}{\rm{op}}
\newcommand{\cop}{\rm{cop}}

\newcommand{\diag}{\,\text{diag}}

\newcommand{\ku}{{\Bbbk}}

\newcommand{\Z}{{\mathbb Z}}
\newcommand{\Na}{{\mathbb N}}
\newcommand{\Di}{{\mathbb D}}

\newcommand{\uno}{{\mathbf 1}}
\newcommand{\C}{{\mathbb C}}

\newcommand{\id}{\mbox{\rm id\,}}

\newcommand{\brp}{\mbox{\rm BrPic\,}}
\newcommand{\Aut}{\mbox{\rm Aut\,}}

\newcommand{\Fun}{\operatorname{Fun}}

\newcommand\Rep{\operatorname{Rep}}

\newcommand\co{\operatorname{co}}

\newcommand{\Id}{\mathop{\rm Id}}

\newcommand{\gr}{\mbox{\rm gr\,}}

\renewcommand{\_}[1]{\mbox{$_{\left( #1 \right)}$}}

\theoremstyle{plain}
\numberwithin{equation}{section}

\newtheorem{teo}{Theorem}[section]
\newtheorem{lema}[teo]{Lemma}

\newtheorem{prop}[teo]{Proposition}
\newtheorem{claim}{Claim}[section]

\theoremstyle{definition}
\newtheorem{defi}[teo]{Definition}
  \newtheorem{exa}[teo]{Example}
\theoremstyle{remark}

\newtheorem{rmk}[teo]{Remark}

\def\pf{\begin{proof}}
\def\epf{\end{proof}}
\theoremstyle{remark}

\begin{document}

\title[The Brauer-Picard group of finite supergroup algebras ]
{The Brauer-Picard group of the representation category of finite supergroup algebras}
\author[ Mombelli]{
Mart\'\i n Mombelli }
\address{ 
Fachbereich Mathematik, Universit\"{a}t Hamburg \newline\indent
Bereich Algebra und Zahlentheorie \newline \indent
Bundesstrasse 55, D –- 20146 Hamburg,
Germany}
\email{martin10090@gmail.com, mombelli@mate.uncor.edu
\newline \indent\emph{URL:}\/ http://www.mate.uncor.edu/$\sim$mombelli}

\begin{abstract} We develop  further the techniques presented
in \cite{M3} to study  bimodule categories over  the representation
categories of arbitrary finite-dimensional Hopf algebras. 
We compute
 the Brauer-Picard group of
equivalence classes of  exact invertible bimodule categories
over  the representation categories of a certain large family of  pointed non-semisimple Hopf algebras, the so called
\emph{supergroup algebras} \cite{AEG}. To obtain this result we first give a
classification of equivalence classes
of exact indecomposable bimodule categories over  such Hopf algebras.

\bigbreak
\bigbreak
\bigbreak
\bigbreak

\end{abstract}

\date{\today}
\maketitle
\setcounter{tocdepth}{1} \tableofcontents
\section{Introduction}

The Brauer-Picard group $\text{BrPic}(\ca)$ of a finite tensor category $\ca$ introduced in
  \cite{ENO},  is the group of equivalence classes 
of invertible exact $\ca$-bimodule categories. This group
is a fundamental piece of information needed to compute extensions of
a given tensor category by a finite group. Also it has a close relation to certain  structures
appearing in mathematical physics, see for example \cite{DKR}, \cite{KK}.

\medbreak

In  \cite{ENO} the authors compute the Brauer-Picard group of
the representation category of an arbitrary finite Abelian group $G$. Given two
semisimple bimodule categories $\Mo, \No$ over $\Rep(G)$
the authors compute the decomposition into indecomposable bimodule categories
of the tensor product $\Mo \boxtimes_{\Rep(G)}\No$. Using
this result and some other techniques they compute 
$\text{BrPic}(\Rep(G))$. The same methods appear to be
unsuccessful for an arbitrary finite-dimensional Hopf algebra $H$. The problem
of explicitly given a decomposition of the tensor 
product $\Mo \boxtimes_{\Rep(H)}\No$ into indecomposable bimodule categories
for arbitrary bimodule categories $\Mo$, $\No$
looks complicated.

\medbreak

Using Hopf theoretic techniques this problem was partially solved in \cite{M3}
by considering the tensor product $\Mo \boxtimes_{\Rep(H)}\No$ only
in the case when both bimodule categories $\Mo$, $\No$ are invertible.

\medbreak

The main result of this paper is the computation of the 
Brauer-Picard group of
the representation category of the so called \emph{supergroup algebras}. 

\medbreak

Let  $G$ be a finite  group, $u$ be an element of order 2 in the center of $G$
and $V$ be a finite-dimensional $G$-module such that $u$ acts by $-1$ in $V$. 
The vector space $V$ is a Yetter-Drinfeld module over $G$ by declaring the coaction
$\delta:V\to \ku G\otk V$, $\delta(v)=u\ot v$, $v\in V$. The Nichols
algebra of $V$ is the exterior algebra $\wedge(V)$ and the bosonization
$\wedge(V)\#\ku G$ is called a \emph{supergroup algebra} \cite{AEG}. We shall denote this
Hopf algebra by $\Ac(V,u,G)$. This family of Hopf algebras played a central role in the
classification of finite-dimensional triangular Hopf algebras \cite{EG}.
\medbreak

If $H$ is a finite-dimensional Hopf algebra then left module categories 
over $\Rep(H)$ are parametrized by equivalence classes of
certain $H$-comodule algebras. Since bimodule categories over $\Rep(H)$
are the same as left module categories over the Deligne's tensor \cite{De}
product $\Rep(H)\boxtimes \Rep(H)^{\op}=\Rep(H\otk H^{\cop}),$ then  
bimodule categories over $\Rep(H)$ are parametrized by equivalence classes of certain
left $H\otk H^{\cop}$-comodule algebras. If $\Mo$ and $\No$ are
invertible exact $\Rep(H)$-bimodule categories the tensor product 
$\Mo \boxtimes_{\Rep(H)}\No$ is an  invertible exact $\Rep(H)$-bimodule category,
therefore indecomposable. In Section \ref{b-tib} we collect all these results and we
recall results from \cite{M3} allowing us to give a precise description
of the category $\Mo \boxtimes_{\Rep(H)}\No$. 

\medbreak

If $H$ is a coradically graded Hopf algebra then $H\otk H^{\cop}$ is also
coradically graded, and indecomposable exact left module categories over $\Rep(H\otk H^{\cop})$
are parametrized by certain equivalence classes of deformations of coideal subalgebras in
$H\otk H^{\cop}$. This results are contained in Section \ref{graded-comod-algs}.

\medbreak

If $\Mo$ is an exact indecomposable bimodule category over $\Rep(\Ac(V,u,G))$
then there exists a certain  left $\Ac(V,u,G)\otk \Ac(V,u,G)^{\cop}$-comodule
algebra $K$ such that $\Mo$ is equivalent to the category of finite-dimensional left $K$-modules.
Since $\Ac(V,u,G)$ is a coradically graded Hopf algebra then $K$ is a certain
deformation of a coideal subalgebra of $\Ac(V,u,G)\otk \Ac(V,u,G)^{\cop}$. 
In Section \ref{defi:qs} we explicitly describe coideal subalgebras in the tensor product
$\Ac(V,u,G)\otk \Ac(V,u,G)^{\cop}$. Using these results, in Section \ref{section:modca}, we prove that if $\Mo$ is an
exact indecomposable left module
category over the category $\Rep(\Ac(V,u,G)\otk \Ac(V,u,G)^{\cop})$ there exists
a 6-tuple $(W^1,W^2,W^3,\beta,F,\psi)$ where
\begin{itemize}
 \item[(i)]  $F\subseteq G\times G$ is a subgroup,   $\psi\in Z^2(F,\ku^{\times})$ is a 2-cocycle,
 \item[(i)]  $W^1, W^2\subseteq V$ $W^3 \subseteq V\oplus V$ are subspaces
 such that $W^3 \cap W^1\oplus W^2= 0$,
$W^3 \cap V\oplus 0=0=W^3 \cap 0\oplus V$, and all subspaces are invariant under the action
of $F$,
 \item[(ii)]  $\beta: \oplus_{i=1}^3 W^i\times \oplus_{i=1}^3 W^i \to \ku$ is a bilinear form
stable under the action of $F$, such that
$$\beta(w_1,w_2)=-\beta(w_2,w_1),\; \beta(w_1,w_3)=\beta(w_3,w_1),\; \beta(w_2,w_3)=-\beta(w_3,w_2), 
$$
for all $w_i\in W^i$, $i=1,2,3$, and $\beta$ restricted to $W^i\times W^i$ is symmetric
for any $i=1,2,3$. If $(u,u)\notin F$ then  $\beta$ restricted to $W^1\times W^2$ is null,
\end{itemize}
such that $\Mo$ is module equivalent to the category 
of finite-dimensional left $\Kc(W^1,W^2,W^3,\beta,F,\psi)$-modules, where 
$\Kc(W^1,W^2,W^3,\beta,F,\psi)$ is a certain left comodule algebra over $\Ac(V,u,G)\otk \Ac(V,u,G)^{\cop}$.
We also describe equivalence classes of such module categories.

\medbreak

Using these results, in Section \ref{sectio:brp}, we prove our main result:
\begin{teo} Assume $G$ is Abelian. The group $\text{BrPic}(\Rep(\Ac(V,u,G)))$ is isomorphic to the
  group of (certain equivalence classes of) pairs  $(T, \alpha)$ where
\begin{itemize}
 \item $\alpha\in O(G\oplus \widehat{G})$, see Definition \ref{orth},
 \item $T:V\oplus V^*\to V\oplus V^*$ is a linear isomorphism such that
$$T(v,f)=x^{-1}\cdot T(y\cdot v,y\cdot f), $$
$$ T^1(0,f)=0, \quad T^2(0,f)(T^1(v,0))=f(v), $$
for all $(v,f)\in  V\oplus V^*$, $(x,y)\in U_\alpha$. Here $T(v,f)=(T^1(v,f), T^2(v,f))$ for all $f\in V^*, v\in V$.
\end{itemize} 
The product of two such triples $(T, \alpha)$, $(T', \alpha')$ is
$$(T, \alpha)\bullet (T', \alpha')=(T\circ T', \alpha\alpha').$$
\end{teo}
As expected, this group is not finite, as is the case for fusion categories.
The main difficulty to prove this theorem relies on finding which of the comodule algebras
$\Kc(W^1,W^2,W^3,\beta,F,\psi)$ give invertible bimodule categories and give an explicit
description of the product of the group $\text{BrPic}(\Rep(\Ac(V,u,G)))$.  Most of Section
\ref{sectio:brp} is dedicated to this task.

\medbreak

It is expected that this result led us to construct interesting new families of finite
non-semisimple tensor categories that are extensions by a finite group
of the category $\Rep(\Ac(V,u,G))$.

\subsection*{Acknowledgment} This work  was  supported by
 CONICET, Secyt (UNC), Mincyt (C\'ordoba) Argentina. It was written during a
 research fellowship granted by CONICET, Argentina in the University of Hamburg, Germany.
 I would like to
thank professors Christoph Schweigert and Ingo Runkel for many  conversations
from which I have benefited greatly and thanks also to professor Bojana Femic for her 
many helpful remarks.

\section{Notation and preliminaries}

We shall work over an algebraically closed field $\ku$ of characteristic 0. All
vector spaces and algebras are considered over $\ku$. We denote $vect_\ku$ the
category of finite-dimensional $\ku$-vector spaces. If $A$ is an algebra we shall denote
by ${}_A\Mo$ ($\Mo_A$ ) the category of finite-dimensional left (right) $A$-modules.

If $V$ is a vector space any bilinear form $\beta:V\times V\to \ku$ determines a linear morphism
$\widehat{\beta}:V \to V^*$
\begin{equation}\label{bilinear-f1}\widehat{\beta}(v)(w)=\beta(v,w),\quad \text{ for all } v,w\in V.
\end{equation}
Let $\Mo$ be an Abelian category. A full subcategory  $\No$ of $\Mo$
is called a \emph{Serre subcategory} if
\begin{itemize}
 \item every object in $\Mo$ isomorphic to an object in $\No$ is in $\No$,
\item every $\Mo$-quotient and every $\Mo$-subobject of an object in $\No$ lies in $\No$,
\item every $\Mo$-extension of objects in $\No$ lies in $\No$.
\end{itemize}
It is well-known that if $F:\Mo\to \Mo$ is an exact functor then the full subcategory of objects $N\in \Mo$
such that $F(N)=0$ is a Serre subcategory. This fact will be used without further mention.

\subsection{Finite tensor categories} A \emph{ tensor category over} $\ku$ is a $\ku$-linear Abelian rigid monoidal category.
Hereafter all tensor categories will be assumed to be over a field $\ku$. A
\emph{finite  category}  is an Abelian $\ku$-linear category such that it has only a finite number of isomorphism
classes of simple objects, Hom spaces are finite-dimensional
$\ku$-vector spaces, all objects have finite lenght and every simple object has a projective cover.
A  \emph{finite tensor category} \cite{eo} is a tensor category with finite underlying Abelian
category such that the unit object is simple. All  functors will be assumed to be
$\ku$-linear and all categories will be finite.

\subsection{Twisting comodule algebras}

Let $H$ be a Hopf algebra. Let us recall that a Hopf 2-cocycle for $H$ is a  map $\sigma: H\otk H\to
\ku$, invertible with respect to convolution,  such that
\begin{align}\label{2-cocycle}
\sigma(x\_1, y\_1)\sigma(x\_2y\_2, z) &= \sigma(y\_1,
z\_1)\sigma(x, y\_2z\_2),
\\
\label{2-cocycle-unitario} \sigma(x, 1) &= \varepsilon(x) =
\sigma(1, x),
\end{align}
for all $x,y, z\in H$. Using this cocycle there is a new Hopf
algebra structure constructed over the same coalgebra $H$ with the
product described by
\begin{equation}\label{product-twisted}
 x._{[\sigma]}y = \sigma(x_{(1)}, y_{(1)}) \sigma^{-1}(x_{(3)},
y_{(3)})\, \, x_{(2)}y_{(2)}, \qquad x,y\in H.
\end{equation}
This new Hopf algebra is denoted by  $H^{[\sigma]}$.
 If $\sigma: H\otimes H\to \ku$ is a Hopf 2-cocycle and $A$ is a left
$H$-comodule algebra, then we can define a new product in $A$ by
\begin{align}\label{sigma-product} a._{\sigma}b = \sigma(a_{(-1)},
b_{(-1)})\, a_{(0)}.b_{(0)},
\end{align}
$a,b\in A$. We shall denote by $A_{\sigma}$ this new algebra. 

\begin{lema} The algebra $A_{\sigma}$ is a left $H^{[\sigma]}$-comodule
algebra.\qed
\end{lema}

\section{Representations of finite tensor categories}\label{section-rep}
Let $\ca$ be a tensor category. For the definition and basic notions of left and right exact module categories
we refer to \cite{eo, O1}.

\medbreak

In this paper we only consider module categories that are finite categories. A module functor between left $\ca$-module categories $\Mo$ and $\Mo'$ over a
tensor category $\ca$ is a pair $(T,c)$, where $T:\Mo \to
\Mo'$ is a  functor and $c_{X,M}: T(X\otb M)\to
X\otb T(M)$ is a family of natural isomorphism such that for any $X, Y\in
\ca$, $M\in \Mo$:
\begin{align}\label{modfunctor1}
(\id_X\otimes c_{Y,M})c_{X,Y\otb M}T(m_{X,Y,M}) &=
m_{X,Y,T(M)}\, c_{X\otimes Y,M}
\\\label{modfunctor2}
\ell_{T(M)} \,c_{\uno ,M} &=T(\ell_{M}).
\end{align}

The  direct sum of two module categories $\Mo_1$ and $\Mo_2$ over
a tensor category $\ca$  is the $\ku$-linear category $\Mo_1\times
\Mo_2$ with coordinate-wise module structure. A module category is
{\em indecomposable} if it is not equivalent to a direct sum of
two non trivial module categories. Any exact module category is 
equivalent to a direct sum of indecomposable exact module categories, see \cite{eo}.

\begin{defi}\cite{AF, Ga1} Let $\Mo$ be a left $\ca$-module category.  A \emph{submodule category}
of $\Mo$ is a Serre subcategory stable under the action of $\ca$.
\end{defi}

The next Lemma is a straightforward consequence of the definitions.
\begin{lema}\begin{itemize}
             \item[1.]  Let $\Mo$ be an exact $\ca$-module category and $\No\subseteq \Mo$
a submodule category. If $\Mo=\oplus_{i\in I} \Mo_i$ is a decomposition into indecomposable
module categories then there is a subset $J\subseteq I$ such that $\No=
\oplus_{i\in J} \Mo_i$.
\item[2.] If  $\Mo$ is an indecomposable exact 
$\ca$-module category and $(F, c):\No\to  \Mo$ is a $\ca$-module functor such that
$F$ is full and faithful, and the subcategory $F(\No)$ is Serre then $F$ is an equivalence.

            \end{itemize}

\end{lema}\qed

\subsection{Bimodule categories}

Let $\ca, \Do$ be tensor categories. For the definition of a $(\ca, \Do)$-\emph{bimodule category}
we refer to \cite{Gr}, \cite{ENO}. A $(\ca, \Do)$-bimodule category
is the same as left $\ca\boxtimes \Do^{\op}$-module category. Here $\boxtimes$ denotes Deligne's
tensor product of Abelian categories \cite{De}. 

A $(\ca, \Do)$-bimodule category is \emph{decomposable} if it is the direct sum of two non-trivial
$(\ca, \Do)$-bimodule categories. A $(\ca, \Do)$-bimodule category is  \emph{indecomposable} 
if it is not decomposable. A $(\ca, \Do)$-bimodule category is \emph{exact}
if it is exact as a left $\ca\boxtimes \Do^{\op}$-module category.
\medbreak

If $\Mo$ is a   right $\ca$-module category then $\Mo^{\op}$ denotes the
 opposite  Abelian category with left $\ca$ action $\ca\times\Mo^{\op} \to \Mo^{\op}$, $(M, X)\mapsto  M\otb X^*$ and associativity
 isomorphisms $m^{\op}_{X,Y,M}=m^{-1}_{Y^*, X^*, M}$ for all $X, Y\in \ca, M\in \Mo$. Similarly if  $\Mo$ is a
   left $\ca$-module category. If  $\Mo$ is a $(\ca,\Do)$-bimodule category then $\Mo^{\op}$ is a 
$(\Do,\ca)$-bimodule category. See \cite[Prop. 2.15]{Gr}.

\medbreak

If $\Mo,  \No$ are $(\ca, \Do)$-bimodule categories, a \emph{bimodule functor}
is the same as a module functor of $\ca \boxtimes \Do^{\op}$-module categories,
that is a  functor $F:\Mo\to\No$ such that  $(F,c):\Mo\to\No$ is a functor of left
 $\ca$-module categories,  also $(F,d):\Mo\to \No$ is a functor of right
$\Do$-module categories and
\begin{equation}\label{bimod-funct-def}
(\id_X\ot d_{M,Y}) c_{X,M\otb_r Y} F(\gamma_{X,M,Y})=
\gamma_{X,F(M),Y} (c_{X,M}\ot \id_Y) d_{X\otb_l M, Y},
\end{equation}
for all $M\in  \Mo$, $X\in\ca$, $Y\in \Do$.

\medbreak

A $(\ca, \Do)$-bimodule category $\Mo$ is called \emph{invertible} \cite[Prop. 4.2]{ENO} if 
there are  equivalences of bimodule categories 
$$\Mo^{\op}\boxtimes_{\ca} \Mo\simeq \Do, \quad \Mo\boxtimes_{\Do} \Mo^{\op}\simeq \ca.$$

\begin{lema}\cite[Corollary 4.4]{ENO} If $\Mo$ is an invertible $(\ca, \Do)$-bimodule category then it is indecomposable
as a bimodule category.\qed
\end{lema}

\begin{lema} \cite[Prop 4.2]{ENO} Let $\Mo$ be an exact  $(\ca, \Do)$-bimodule category. The following
statements are equivalent.
\begin{itemize} 
\item[1.] $\Mo$ is an invertible.
\item[2.]  There exists a $\Do$-bimodule equivalence $\Mo^{\op}\boxtimes_{\ca} \Mo\simeq \Do$.
\item[3.] There exists a $\ca$-bimodule equivalence $ \Mo\boxtimes_{\Do} \Mo^{\op}\simeq \ca$.
\item[4.] The functor $R:\Do^{\op}\to \Fun_\ca(\Mo,\Mo)$, $R(X)(M)= M\otb X$, for all
$X\in\Do$, $M\in \Mo$, is an equivalence of tensor categories.
\item[2.] The functor $L:\ca\to \Fun_\Do(\Mo,\Mo)$, $L(Y)(M)=Y\otb M$, for all
$Y\in\ca$, $M\in \Mo$, is an equivalence of tensor categories.
            \end{itemize}
\end{lema}
\pf The proof of \cite[Prop 4.2]{ENO} extends \emph{mutatis mutandis} to the non-semisimple
case using results from \cite{eo}.
\epf

\subsection{Module categories over Hopf algebras}\label{modcat-over-hopf}
 
Let $H$ be a finite-dimensional Hopf algebra
 and let $(\Ac,\lambda)$ be a left $H$-comodule algebra.  The
category  ${}_\Ac\Mo$ is a representation
of $\Rep(H)$. The action 
$$\otb:\Rep(H)\times {}_\Ac\Mo\to {}_\Ac\Mo, \quad V\otb M=V\otk M,$$ 
for all $V\in \Rep(H)$, $M\in {}_\Ac\Mo$. The left $\Ac$-module
structure on $V\otk M$ is given by 
$$a\cdot (v\ot m)=a\_{-1}\cdot v\ot a\_0\cdot m,$$
for all $a\in \Ac$, $v\in V$, $m\in M$. Here $ \lambda:\Ac\to H\otk \Ac$, $\lambda(a)=a\_{-1}\ot a\_0$.
\medbreak
If $\Ac$ is a $H$-comodule algebra via $\lambda:\Ac\to H\otk \Ac$, we
shall say that a (right) ideal $J$ is $H$-costable if
$\lambda(J)\subseteq H\otk J$. We shall say that $\Ac$ is (right)
$H$-simple, if there is no nontrivial (right) ideal $H$-costable
in $\Ac$. When $\Ac$ is right $H$-simple then the category ${}_\Ac\Mo$ is exact. 

\begin{teo}\cite[Theorem 3.3]{AM}\label{mod-over-hopf} Let $\Mo$ be an exact indecomposable module category over
 $\Rep(H)$ then
there exists a left $H$-comodule algebra $\Ac$ right $H$-simple with trivial
coinvariants such that $\Mo\simeq {}_\Ac\Mo$ as $\Rep(H)$-modules.

\qed
\end{teo}

Two $H$-comodule algebras $\Ac$, $\Ac'$ are \emph{equivariantly Morita equivalent} if 
the module categories  ${}_{\Ac'}$, ${}_\Ac\Mo$ are equivalent.

\section{Bimodule categories over Hopf algebras}\label{b-tib}

\subsection{Tensor product of invertible bimodule categories}\label{tib}

Let $A, B$ be  finite-dimensional Hopf algebras. A $(\Rep(B), \Rep(A))$-bimo\-dule category is the same as a left
$\Rep( B\otk A^{\cop})$-module category. This follows from the fact
that $\Rep(A)^{\op}\simeq \Rep(A^{\cop})$ and $ \Rep(B)\boxtimes   \Rep(A^{\cop})\simeq
\Rep(B\otk A^{\cop})$. Thus  Theorem \ref{mod-over-hopf} implies that
 any exact indecomposable $(\Rep(B), \Rep(A))$-bimodule category is equivalent
to the category ${}_S\Mo$ of finite-dimensional left $S$-modules, where $S$ is a finite-dimensional  right $B\otk A^{\cop}$-simple left $B\otk A^{\cop}$-comodule
algebra. 

\medbreak

Since $\Rep(A)$ is canonically a $\Rep(A)$-bimodule category then there exists
some right $A\otk A^{\cop}$-simple left $A\otk A^{\cop}$-comodule
algebra $\Ac$ such that $\Rep(A)\simeq {}_\Ac\Mo$ as $\Rep(A)$-bimodule categories.
In \cite{M3} we computed this comodule algebra. Let us recall this result.

We denote by  $\diag(A)$ the left
$A\otk A^{\cop}$-comodule algebra with  underlying algebra $A$ and  comodule
structure:
$$\lambda:\diag(A)\to A\otk A^{\cop}\otk \diag(A), \quad \lambda(a)=a\_1\ot a\_3\ot a\_2,$$
for all $a\in A$.
Thus the category ${}_A\Mo$ is a $\Rep(A)$-bimodule category.
\begin{lema}\label{diagonal-comod} \begin{itemize}
                                    \item[1.] $\diag(A)$ is a right simple left $A\otk A^{\cop}$-comodule
algebra and $\diag(A)^{\co A\otk A^{\cop}}=\ku 1.$
 \item[2.]  There is an equivalence of $\Rep(A)$-bimodule categories
$$ {}_A\Mo\simeq {}_{\diag(A)}\Mo.$$
                                   \end{itemize}
\end{lema}
\pf 1. Let $0\neq I\subseteq A$ be a right ideal $A$-costable. Then for any $a\in I$, $a\_1\ot a\_3\ot a\_2
\in A\otk A\ot I$ which implies that $a\_1\ot a\_2\in A\otk I$. Thus $I$ is a right ideal 
stable under the coaction, then $I=A$. 

2. The identity functor $\Id: {}_A\Mo\to  {}_{\diag(A)}\Mo$ is an equivalence of  $\Rep(A)$-bimodule categories.
\epf

Let us recall some constructions and results  obtained in \cite{M3} concerning the tensor product of bimodule
categories over Hopf algebras.  Set
$\pi_A:A\ot B\to A$, $\pi_B:A\ot B\to B$  the algebra maps defined by
$$\pi_A(x\ot y)=\epsilon(y) x, \quad \pi_B(x\ot y)=\epsilon(x) y,$$
for all $x\in A, y\in B.$ 
\medbreak

Let $K$ be a right $B\ot A^{\cop}$-simple left $B\ot A^{\cop}$-comodule algebra and
$L$ a right $A\ot B^{\cop}$-simple left $A\ot B^{\cop}$-comodule algebra. Thus
the category ${}_K\Mo$  is a $(\Rep(B), \Rep(A))$-bimodule category and
  ${}_L\Mo$ is a  $(\Rep(A), \Rep(B))$-bimodule category. 

The category $\Mo(A,B,K,\overline{L})$ is the category ${}_{K}^B\Mo_ {\overline{L}}$ 
of $(K,\overline{L})$-bimodules and left $B$-comodules such that the comodule
structure is a bimodule morphism. See
 \cite[Section 3]{M3}. It has a structure of $(\Rep(A), \Rep(A))$-bimodule category.
Recall that $\overline{L}$ is the  left $B\ot A^{\cop}$-comodule algebra
with opposite algebra structure  $L^{\op}$ and  left $B\ot A^{\cop}$-comodule
structure:
\begin{equation}\label{comod-op} \overline{\lambda}:L\to A^{\cop}\otk B\otk L,\quad
 l\mapsto  (\Ss^{-1}_B\ot \Ss_A)(l\_{-1})\ot l\_0,
\end{equation}
for all $l\in L$.
Also $L$ is a right $B$-comodule with comodule map given by
\begin{equation}\label{comod-op2} l\mapsto l\_0\ot \pi_B(l\_{-1}),
\end{equation}
 for all $l\in L$, and $K$ is a left $B$-comodule 
with comodule map given by
\begin{equation}\label{comod-op3} k\mapsto  \pi_B(k\_{-1})\ot k\_0,
\end{equation}
for all $k\in K$. Using this structure we can form the cotensor product $L
\Box_B K$. Define 
\begin{equation}\label{coact-coprodu-t} \lambda(l\ot k)=\pi_A(l\_{-1})\ot \pi_A(k\_{-1}) \ot l\_0\ot k\_0,
\end{equation}
for all $l\ot k\in L \Box_B K$. Then
$L \Box_B K$ is a left $A\otk A^{\cop}$-comodule algebra. See \cite[Lemma 3.6]{M3}.

\medbreak

In  \cite{M3} we have presented the functors
$$\Fc: {}_{L \Box_B K}\Mo\to \Mo(A,B,K,\overline{L}),\quad \Gc:\Mo(A,B,K,\overline{L})\to {}_{L \Box_B K}\Mo$$
 by $\Fc(N)=(L\otk K)\ot_{L \Box_B K}  N$ 
for all $N\in {}_{L \Box_B K}\Mo$ and  $\Gc(M)=M^{\co B}$ for all
$M\in \Mo(A,B,K,\overline{L})$. We recall that the left $B$-comodule structure on
$\Fc(N)$ is given by $\delta: \Fc(N)\to B\otk \Fc(N)$,
\begin{equation}\label{comod-str-B} 
 \delta(l\ot k\ot n)=\pi_B(k\_{-1})\Ss^{-1}(\pi_B(l\_{-1}))\ot l\_0\ot k\_0\ot n,
\end{equation}
for all $l\in L, k\in K, n\in N.$

\medbreak

This pair of functors were studied in  \cite{D}, \cite{CCMT}. In the following theorem  we summarize  some results from \cite{M3}.
\begin{teo}\label{tensor-bimod-hopf}  \begin{itemize}
  \item[(a)] There is a  $\Rep(A)$-bimodule equivalence:
$$ {}_L\Mo \boxtimes_{\Rep(B)}
 {}_K\Mo \simeq  \Mo(A,B,K,\overline{L}).$$
\item[(b)] $\Fc$ and $\Gc$ are   $\Rep(A)$-bimodule functors.

\item[(c)] Assume that both bimodule categories $ {}_L\Mo$, ${}_K\Mo$ are invertible and
$L\otk K\simeq C\otk  L \Box_B K$, as  right  $L \Box_B K$-modules and  left $B$-comodules. Here
$C$ is a certain left $B$-comodule.
Then there is
an equivalence of $\Rep(A)$-bimodule categories
$${}_{L \Box_B K}\Mo \simeq  {}_L\Mo \boxtimes_{\Rep(B)}
 {}_K\Mo.$$

\end{itemize}
 
\end{teo}
\pf For the proof of (a) and (b)  see \cite{M3}.

\medbreak

(c). We shall prove that the functors $\Fc$, $\Gc$ establish an equivalence
 of module categories.

Let us prove  that $\Fc(\Gc(M))\simeq M$ for all $M\in  \Mo(A,B,K,\overline{L})$.
 For any $M\in  \Mo(A,B,K,\overline{L})$ there is a projection
$$ \pi_M: (L\otk K)\ot_{L \Box_B K} M^{\co B} \to M, \quad \pi_M(l\ot k\ot m)=(l\ot k)\cdot m,$$
for all $l\ot k\in L\otk K$, $m\in M^{\co B}$. Define the functor $\Phi: \Mo(A,B,K,\overline{L})
\to vect_\ku$, $\Phi(M)=\ker(\pi_M)$. The functor $\Phi$ is  a module functor. To see this it is enough
to prove that the diagram
\begin{align*}
\begin{CD}
\Fc(\Gc(X\otb M)) @>\simeq >> X\otb \Fc(\Gc( M)) \\
@V\pi_{X\otb M}VV    @VV\id_X\ot \pi_M V  \\
 X\ot M @>>\id> X\ot M.
\end{CD}
\end{align*}
is commutative. Then $\Phi$ is exact. The full subcategory $\No$ of $\Mo(A,B,K,\overline{L})$ consisting of objects $M$
such that $\Phi(M)=0$ is a submodule category. $\No$ is not the null category
since $\pi_{L \otk K}=\id$, thus $L \otk K\in\No$. Since both $ {}_L\Mo$, ${}_K\Mo$ are invertible
the product $ {}_L\Mo \boxtimes_{\Rep(B)}
 {}_K\Mo \simeq  \Mo(A,B,K,\overline{L})$ is indecomposable. Hence $\No
= \Mo(A,B,K,\overline{L}).$ This implies that $\Fc(\Gc(M))=M$ for all $M\in  \Mo(A,B,K,\overline{L})$.
Since $L\otk K\simeq C\otk  L \Box_B K$, as  right  $L \Box_B K$-modules and  left $B$-comodules
the functor $\Fc$ is full and faithful, thus it is an equivalence of categories.\epf

\begin{rmk} In all examples the assumption $L\otk K\simeq C\otk  L \Box_B K$
in Theorem \ref{tensor-bimod-hopf} (c) seems to be
superfluous, although I do not know any counterexample. 
\end{rmk}

\section{Graded comodule algebras over Hopf algebras}\label{graded-comod-algs}

From the discussion on Section \ref{modcat-over-hopf} equivalence classes
 of indecomposable exact module categories over the representation categories of Hopf algebras are
in correspondence with equivariant Morita equivalence classes of right simple comodule algebras.
To study this class of algebras we developed a technique  in \cite{M1} 
using the Loewy filtration and the associated graded algebra. We briefly recall all
this notions.

\medbreak

If $H$ is a finite-dimensional Hopf algebra then $H_0\subseteq H_1
\subseteq \dots \subseteq H_m=H$ will denote the coradical
filtration. When $H_0\subseteq H$ is a Hopf subalgebra then the
associated graded algebra $\gr H$ is a coradically graded Hopf
algebra. If $(A, \lambda)$ is a left $H$-comodule algebra, the
coradical filtration on $H$ induces a filtration on $A$, given by
$A_n=\lambda^{-1}(H_n\otk A)$, $n=1,\dots, m$. This filtration is called the
\emph{Loewy series }on $A$.

\medbreak

Recall that if $H=\oplus^m_{i=0} H(i)$ is a coradically graded Hopf algebra, 
a left $H$-comodule algebra $(A,\lambda)$ is a \emph{graded comodule algebra},
 if it is graded as an algebra
$A=\oplus^m_{i=0} A(i)$ and for
each $0\leq n\leq m$
\begin{equation}\label{graded-comod-st}
 \lambda(A(n))\subseteq \bigoplus^m_{i=0} H(i)\otk A(n-i).
\end{equation}
A graded comodule algebra $A=\oplus^m_{i=0} A(i)$ is
\emph{Loewy-graded} if the Loewy series is given by
$A_n=\oplus^n_{i=0} A(i)$ for any $0\leq n\leq m$.

\medbreak

If $A$ is a left $H$-comodule algebra the associated  graded algebra $\gr A$
obtained from the Loewy filtration is a Loewy-graded left $\gr
H$-comodule algebra. For more details see \cite{M1}.

\medbreak

The following result will be needed later. 
\begin{lema}\label{comod-st-tech} Let $H=\oplus^m_{i=0} H(i)$
be a coradically graded Hopf algebra and $(A,\lambda_A)$ a left $H$ comodule with a 
grading $A=\oplus_{i=0}^m A(i)$ such that \eqref{graded-comod-st} holds. If $B\subseteq A$
is a subcomodule algebra and we set $B(n)=B\cap A(n)$ then
$$ B=\oplus_{i=0}^m B(i).$$
\end{lema}
\pf Let $b\in B$, then $b=\sum_{i=0}^m b_i$ where $b_i\in A(i)$. Let us prove that
$b_i\in B$ for all $i=0,\dots, m$. Denote $p:H\to H(0)$, $\pi_j:A\to A(j)$ the
canonical projections. Observe that for any $j=0,\dots, m$
$$(p\ot \pi_j)\lambda(b)= (p\ot \id)\lambda(b_j).$$
Since $(\epsilon\ot \id)(p\ot \id)\lambda(b_j)=b_j$ then 
$b_j=(\epsilon\ot \id)(p\ot \pi_j)\lambda(b)\in B$.
\epf

\subsection{Comodule algebras over coradically graded Hopf algebras}

Let $G$ be a finite group and $H=\oplus_{i=0}^m H(i)$ be a finite-dimensional coradically graded Hopf algebra 
where $H(0)=\ku G$ is the coradical. 

\medbreak

Let $(A,\lambda)$ be a left $H$-comodule algebra
right $H$-simple with trivial coinvariants and with a grading $A=\oplus_{i=0}^m A(i)$
 making $A$ a Loewy-graded
left $H$-comodule algebra. Since $A$ is right $H$-simple with trivial coinvariants then
$A(0)=\ku_{\psi} F$ where $F\subseteq G$ is a subgroup and $\psi\in Z^2(F,\ku^\times)$ is a 
2-cocycle.

\medbreak

Set $\pi:A\to A(0)$ the canonical projection and $\epsilon:A(0)\to \ku$ the map
given by $\epsilon(e_f)=1$ for all $f \in F$. 
\begin{rmk} If $\psi$ is trivial then $\epsilon:A(0)\to \ku$ is an algebra morphism.
\end{rmk}

\begin{prop}\label{comod=coideal} Assume that $\psi$ is trivial and let  $\phi:A\to H$ be the map defined by $\phi=(\id_H\ot \epsilon \pi)
\lambda$. Then
\begin{itemize}
 \item[(i)] $\phi$ is an algebra morphism,
  \item[(ii)] $\phi$ is a $H$-comodule map, and
 \item[(iii)] $\phi$ is injective.
\end{itemize}

\end{prop}
\pf (i).  It follows since all maps in the definition of $\phi$ are algebra morphisms.

(ii). $\Delta  \phi= \Delta (\id_H\ot \epsilon \pi)
\lambda =  (\id_H\ot \id_H\ot\epsilon \pi) (\Delta\ot \id_A) \lambda$. Using the coassociativity
of $\lambda$ we obtain that $\Delta  \phi= (\id_H\ot \phi)\lambda$.

(iii). Let $a\in \ker \phi$. Assume that $a\neq 0$. Write $a=\sum_{n=0}^t a^{(n)}$ where $a^{(n)}\in A(n)$ and $t\leq m$. 
We can assume that  $ a^{(t)}\neq 0$. Then $\lambda(a^{(n)})\in \oplus_{i=0}^n H(i)\otk A(n-i)$.
Set $\lambda(a^{(n)})=\sum_{i=0}^n b_{n,i} $ where $b_{n,i}= \sum_k x_k^{n,i}\ot c_k^{n,i}$
and $x^{n,i}\in H(i),$ $c^{n,i}\in A(n-i)$.

Since $ a^{(t)}\neq 0$ then $b_{t,t}\neq 0$. Indeed, if $b_{t,t}=0$ then
$\lambda(a^{(t)} )\in  \oplus_{i=0}^{t-1} H(i)\otk A $, hence $a^{(t)}\in \oplus_{i=0}^{t-1} 
A(i)=A_{t-1}$, which is imposible unless $ a^{(t)}=0$.

Also,  $\Delta\phi(a)=0$, then $(\id_H\ot \id_H\ot \epsilon \pi) (\id_H\ot \lambda) \lambda(a)=0$
which implies that
\begin{align*}  \sum_{n=0}^t \sum_{i=0}^n  \sum_k  x_k^{n,i}\ot (\id_H\ot \epsilon \pi)\lambda(c_k^{n,i})
=0
\end{align*}
The element of the above summation that belongs to $H(t)\otk H(0)\otk A(0)$ must be equal to
zero, hence
$\sum_k  x_k^{t,t}\ot (\id_H\ot \epsilon)\lambda(c_k^{t,t})=0.$
Since we have that $\sum_k  x_k^{t,t}\ot (\id_H\ot \epsilon)\lambda(c_k^{t,t})= b_{t,t}$ we get
that $b_{t,t}=0$ which is a contradiction.
Therefore $a=0$. 
\epf
In another words, Proposition \ref{comod=coideal} implies that if $A$ is a Loewy-graded right $H$-simple
left comodule algebra with trivial coinvariants and $A(0)$ is a Hopf subalgebra
of $H(0)$ then $A$ is isomorphic to a left coideal subalgebra of $H$. The next
step is to study what happens if $A(0)$ is not a Hopf subalgebra of $H(0)$.
\medbreak
Let $\widehat{\psi}\in Z^2(G,\ku^{\times})$ be a 2-cocycle such that 
$\widehat{\psi}\mid_{F\times F}=\psi$.
\begin{lema}\label{lift-twist} There exists a Hopf 2-cocycle
$\sigma:H\otk H\to \ku$ such that for any homogeneous elements $x, y\in H$
\begin{align}\label{lift-twist-def}
\sigma(x,y)=\begin{cases}
               \widehat{\psi}(x,y), & \text{if } x,y\in H(0);\\
0, &\text{otherwise.}
               \end{cases}\end{align}
\end{lema}
\pf See \cite[Lemma 4.1]{GM}. \epf

The following result is a straightforward consequence
of Proposition \ref{comod=coideal}.
\begin{lema}\label{transp-comod}  Let $A$ be a Loewy-graded right $H$-simple
left comodule algebra with trivial coinvariants and $A(0)=\ku_{\psi} F$ where $F\subseteq G$ is a subgroup and $\psi\in Z^2(F,\ku^\times)$ is a 
2-cocycle. Then, there exists a Hopf 2-cocycle
$\sigma:H\otk H\to \ku$ such that $A_\sigma$ is isomorphic to a homogeneous left
coideal subalgebra of $H^{[\sigma]}$ as a left $H^{[\sigma]}$-comodule algebras.\qed
\end{lema}
\pf From Lemma \ref{lift-twist} there exists a Hopf 2-cocycle 
$\sigma:H\otk H\to \ku$ such that $\sigma(x,y)=\psi^{-1}(x,y)$ for all $x,y\in F$. The comodule
algebra $A_\sigma$ is Loewy-graded and $(A_\sigma)(0)=\ku F$. Thus the Lemma follows
from Proposition \ref{comod=coideal}. \epf

\section{Supergroup algebras and its coideal subalgebras}\label{defi:qs}

We shall recall the definition of supergroup algebras \cite{AEG}, its Hopf
algebra structure and we describe the
tensor product of two such Hopf algebras. We compute also its homogeneous coideal
subalgebras, a key ingredient to compute module categories.

\subsection{Finite supergroup algebras}\label{spa} 
Let $G$ be a finite group, $u\in G$ be a central element of order 2 and $V$ a finite-dimensional 
$G$-module such that $u\cdot v=-v$  for all $v\in V$.  
The space $V$ has a $G$-comodule structure $\delta:V\to \ku G\otk V$ given by
$\delta(v)=u\ot v$, for all $v\in V$. This gives $V$ structure of Yetter-Drinfeld module
over $\ku G$. The Nichols algebra of $V$ is the exterior algebra $\nic(V)=\wedge(V)$.
The Hopf algebra obtained by bosonization $\wedge(V)\# \ku G$ is called in \cite{AEG} a \emph{finite supergroup algebra}.
 We will denote this Hopf algebra by $\Ac(V,u,G)$. Hereafter we shall denote
the element $v\# g$ simply by $vg$, for all $v\in V, g\in G$.

The algebra $\Ac(V,u,G)$ is generated by elements $v\in V, g\in G$ subject to relations
$$vw+wv=0, \quad gv= (g\cdot v) g, \text{ for all } v,w\in V, g\in G.$$
The coproduct and antipode are determined by
$$\Delta(v)=v\ot 1+ u\ot v, \quad \Delta(g)=g\ot g,$$
$$\Ss(v)=-u v, \quad \Ss(g)=g^{-1},$$
 for all $v\in V, g\in G$.

\begin{lema}\label{iso-cop} There is a Hopf algebra isomorphism
$$\Ac(V,u,G)\simeq \Ac(V,u,G)^{\cop}.$$
\end{lema}
\pf Let $\phi: \Ac(V,u,G)\to \Ac(V,u,G)$ be the algebra map determined by
$$\phi(v)=vu, \quad \phi(g)=g,$$
for all $v\in V, g\in G$. It follows by a direct computation that 
$\phi$ is a Hopf algebra isomorphism between $\Ac(V,u,G)$ and $\Ac(V,u,G)^{\cop}.$
\epf

\subsection{Tensor product of supergroup algebras}\label{tp-spa} Let  $G_1, G_2$ be  finite groups and 
$u_i\in G_i$ be central elements of order 2. For $i=1,2$ let $V_i$ be finite-dimensional
$G_i$-modules,  such that $u_i$ acts in $V_i$ as $-1$. We shall describe the tensor product Hopf algebra
$\Ac(V_1,u_1,G_1)\otk \Ac(V_2,u_2,G_2).$
From now on, we shall denote this Hopf algebra by $\Ac(V_1,V_2,u_1,u_2,G_1,G_2)$.
Let us give a presentation by generators and relations of this algebra.

Set $G=G_1\times G_2$. Both vector spaces $V_1, V_2$ are $G$-modules by setting
$$ (g,h)\cdot v_1= g\cdot v_1, \quad   (g,h)\cdot v_2= h\cdot v_2,$$
for all $(g,h)\in G$, $v_i\in V_i$, $i=1,2$.
The algebra $\Ac(V_1,V_2,u_1,u_2,G_1,G_2)$
is generated by elements
$ V_1, V_2, G$
subject to relations
$$  v_1w_1+ w_1v_1=0,\;\;  v_2w_2+ w_2v_2=0,\; \; v_1v_2=v_2v_1,$$
$$g v_1=  g\cdot v_1 g,\quad   gv_2= g\cdot v_2 g,$$
for all $g\in G$, $v_i\in V_i$, $i=1,2$. The Hopf algebra structure
is determined by
$$ \Delta(v_1)= v_1\ot 1+ (u_1,1)\ot v_1, \;\;\; \Delta(v_2)=v_2\ot 1+ (1,u_2)\ot v_2,$$
$$ \Delta(g_1,g_2)=(g_1,g_2)\ot (g_1,g_2),$$
for all $(g_1,g_2)\in G$, $v_i\in V_i$, $i=1,2$.

\medbreak

We shall define a family of Hopf algebras that are cocycle deformations of
tensor product of supergroup algebras. Let $(V_1,V_2,u_1,u_2,G_1,G_2)$ be a data as
above. Set $V=V_1\oplus V_2$. Define $\Hc(V_1,V_2,u_1,u_2,G_1,G_2)= \wedge(V)\otk \ku G$ with product
determined by
$$vw+wv=0, \quad g v= (g\cdot v) g, \text{ for any }\, v,w\in V_1\oplus V_2, g\in G,$$ 
and coproduct determined by
$$ \Delta(v_1)= v_1\ot 1+ (u_1,1)\ot v_1, \;\;\; \Delta(v_2)=v_2\ot 1+ (1,u_2)\ot v_2,$$
for any $v_i\in V_i$, $i=1,2$.

\begin{prop}\label{twisting-supergroup} Let $H= \Ac(V_1,V_2,u_1,u_2,G_1,G_2)$ and $\sigma:H
\otk H\to \ku$ a Hopf 2-cocyle coming from a 2-cocycle $\psi\in Z^2(G,\ku^{\times})$ as in
Lemma \ref{lift-twist}. Denote $\xi=\psi((u_1,1),(1,u_2))\psi((1,u_2),(u_1,1))^{-1}$. Then
\begin{itemize}
 \item[(i)] if $\xi=1$ we have $H^{[\sigma]}\simeq \Ac(V_1,V_2,u_1,u_2,G_1,G_2)$
 \item[(ii)] if $\xi=-1$ then $H^{[\sigma]}\simeq \Hc(V_1,V_2,u_1,u_2,G_1,G_2)$.
\end{itemize}
\end{prop}
\pf Let $v \in V_1, w\in V_2$ then
$$(\id\ot \Delta)\Delta(v)=v\ot 1\ot 1+(u_1,1)\ot v\ot 1+ (u_1,1)\ot (u_1,1)\ot  v,$$
$$(\id\ot \Delta)\Delta(w)=w\ot 1\ot 1+ (1,u_2)\ot w\ot 1+  (1,u_2)\ot  (1,u_2)\ot  w.$$
Therefore, using \eqref{product-twisted}, it follows that for any $v_1,w_1\in V_1$, 
$v_2,w_2\in V_2$
$$v_1\cdot_{[\sigma]} w_1 + w_1\cdot_{[\sigma]} v_1 =0, \quad 
v_2\cdot_{[\sigma]} w_2 + w_2\cdot_{[\sigma]} v_2 =0,$$
$$v_1\cdot_{[\sigma]} w_2 - \xi \; w_2\cdot_{[\sigma]} v_1 =0. $$
Also for any $g\in G$, $i=1,2$ 
$$ g\cdot_{[\sigma]} v_1 =\psi(g , (u_1,1)) \, g v_1,\quad  v_1\cdot_{[\sigma]} g=  \psi( (u_1,1), g)
\, v_1g, $$
$$ g\cdot_{[\sigma]} v_2 = \psi( (1,u_2), g)\, g v_2,\quad  v_2\cdot_{[\sigma]} g= \psi( (1,u_2), g)
\, v_2g. $$
Hence 
$$ g\cdot_{[\sigma]} v \cdot_{[\sigma]} g^{-1}= gvg^{-1},$$
for any $v\in V$. From these relations, and since the coproduct remains unchanged, we deduce that  if $\xi=1$ then
$H^{[\sigma]}\simeq \Ac(V_1,V_2,u_1,u_2,G_1,G_2)$ and if $\xi=-1$ then $H^{[\sigma]}\simeq \Hc(V_1,V_2,u_1,u_2,G_1,G_2)$.
\epf

\subsection{Homogeneous coideal subalgebras in supergroup algebras}
A \emph{homogeneous left coideal subalgebra} of a coradically graded Hopf
 algebra $H=\oplus_{i=0}^m H(i)$ is a  left coideal subalgebra
$K\subseteq H$ together with an algebra grading $K=\oplus_{i=0}^m K(i)$
 such that $K(i)\subseteq
H(i)$. The main goal of this section is the classification of homogeneous
coideal subalgebras in the tensor product of supergroup algebras.

\medbreak

Let  $(V_1,V_2,u_1,u_2,G_1,G_2)$ be a data as in section \ref{tp-spa}.
 Denote $V=V_1\oplus V_2$
and $u=(u_1,u_2)\in G=G_1\times G_2$. Also set $H=\Ac(V_1,V_2,u_1,u_2,G_1,G_2)$
and $\widetilde{H}=\Hc(V_1,V_2,u_1,u_2,G_1,G_2)$. If $(v_1,v_2)\in V$
 we denote  
$$[(v_1,v_2)]=v_1+ v_2 u\in H(1).$$

\begin{rmk}\label{brkt-vect} For any $(v_1,v_2)\in V$ we have
\begin{equation}\label{cop-bracket} 
[(v_1,v_2)]^2=0, \;\; \Delta([(v_1,v_2)])= v_1\ot 1 +  v_2 u\ot u + (u_1,1)\ot [(v_1,v_2)].
\end{equation}
\end{rmk}

\begin{defi} A \emph{coideal subalgebra data} is  a collection $(W^1,W^2,W^3, F)$, where
\begin{itemize}
 \item $W^1\subseteq V_1$ and $
W^2\subseteq V_2$ are subspaces,
 \item   $W^3\subseteq V$ is a subspace such that $W^3 \cap W^1\oplus W^2= 0$,
$W^3 \cap V_1=0=W^3 \cap V_2$,
 \item $F\subseteq G$ is  a subgroup that leaves invariant all  subspaces $W^i$, $i=1,2,3$,
 \item if $W^3\neq 0$ we require that $u\in F$.
\end{itemize}
we denote $C(W^1,W^2,W^3, F)$ the subalgebra
of $H$ generated by $\ku F$ and elements in $W^1\oplus W^2$ and $\{[w]: w\in W^3\}$. 

\end{defi}

\begin{lema} The algebra $C(W^1,W^2,W^3, F)$ is a homogeneous left
coideal subalgebra of $H$.\qed
\end{lema}

\begin{teo}\label{coideal-tensor-prod-supergroup} Let $K=\oplus_{i=0}^m K(i) \subseteq H$ be a homogeneous left coideal subalgebra. There
exists a coideal subalgebra data $(W^1,W^2,W^3, F)$ such that
$K=C(W^1,W^2,W^3, F)$.
\end{teo}
\pf 
Since $K(0)
\subseteq \ku G$ is a left coideal subalgebra then $K(0)=\ku F$ for some subgroup
$F\subseteq G$. 
 If $K(1)=0$ then $K=\ku F$. Indeed, if $x\in K(2)$ then $\Delta(x)\in H(0)\otk K(2) \oplus H(2)\ot K(0)$,
hence $x\in H_1$ and since $H_1\cap H(2)=0$ then $x=0$. In a similar way we can prove that $K(n)=0$ for all $n$.
 \medbreak
Thus we can assume that $K(1)\neq 0$. The vector space $K(1)$ is a $\ku G$-subcomodule of $V\otk \ku G$ via
$$(\pi\ot \id) \Delta:K(1) \to  \ku G\otk K(1),$$
where $\pi:H\to \ku G$ is the canonical projection. Thus $K(1)=\oplus_{g\in G}
K(1)_g$, where  $K(1)_g=\{k\in K(1): (\pi\ot \id) \Delta(k)=g\ot k\}$, and
$$K(1)_g\subseteq  V_1\otk \ku <(u_1,1)g> \oplus\, V_2\otk \ku <(1,u_2)g>.$$
Therefore we can write 
$$K(1)_{(u_1,1)}= W^1\oplus  \widetilde{W}^2 u \oplus U^3,$$
$$K(1)_{(1,u_2)}= W^2\oplus  \widetilde{W}^1 u \oplus  \widetilde{U}^3,$$  where
$W^1$ is the intesection of $K(1)_{(u_1,1)}$ with $V_1$,  
$\widetilde{W}^2$ is the intesection  of $K(1)_{(u_1,1)}$ with $V_2\otk \ku <u>$ and $U^3$ is a
direct complement, that is, a vector subspace of $V_1 \oplus\, V_2\otk  \ku <u>$ consisting
of elements of the form $[w]$ where $w\in W^3$ and $W^3\subseteq V_1\oplus V_2$.
Since $U^3 \cap W^1\oplus \widetilde{W}^2 u=0$ then $W^3 \cap W^1\oplus \widetilde{W}^2= 0$.
The same is done for $K(1)_{(1,u_2)}$, that is $W^2$ is the intersection
of $K(1)_{(1,u_2)}$ with $V_2$, $ \widetilde{W}^1 u$ is the intersection of $K(1)_{(1,u_2)}$ with
with  $V_1\otk \ku <u>$ and $ \widetilde{U}^3$ is a direct complement. The space
$ \widetilde{U}^3$ consists of elements of the form 
$[w]$ where $w\in \widetilde{W}^3$ and $\widetilde{W}^3\subseteq V_1\oplus V_2$.

\begin{claim} If $u\notin F$ then 
 $ \widetilde{W}^2= \widetilde{W}^1=\widetilde{W}^3=W^3=0$.  We have that $ \widetilde{W}^2= W^2$, 
$\widetilde{W}^1=W^1$.
\end{claim} 
\pf[Proof of Claim] Let $0\neq (v,w)\in W^3$, then $0\neq [(v,w)]\in U^3$. Since
$\Delta( [(v,w)])\in H(0)\otk K(1) \oplus H(1)\otk K(0)$, using \eqref{cop-bracket},
we get that $u\in F$. The same argument works if $\widetilde{W}^3\neq 0$,
 $\widetilde{W}^1\neq 0$ or if  $\widetilde{W}^2\neq 0$.

\medbreak
Let $0\neq w\in \widetilde{W}^2$, then
$wu\in K(1)_{(u_1,1)}$. Since $u\in F$ then $w\in K(1)_{(1,u_2)}$ and the only possibility is that
$w\in W^2$. The other inclusion is  proven similarly. Thus $ \widetilde{W}^2= W^2$. The equality
$\widetilde{W}^1=W^1$ follows analogously. 
\epf

We claim that $K(1)= W^1 F\oplus W^2 F \oplus U^3 F$. Indeed, take $g\in G$ and $0\neq w\in K(1)_g$,
then
$$w= w_1 (u_1,1)g + w_2  (1,u_2)g,$$
for some $w_1\in V_1, w_2\in V_2$. Note that
\begin{align}\label{cop-w}\Delta(w)=w_1  (u_1,1)g\ot (u_1,1)g +  g\ot w +w_2 (1,u_2)g\ot 
 (1,u_2)g.
\end{align}
If $w_1\neq 0$, since  $\Delta(w)\in H(0)\otk K(1) \oplus H(1)\otk K(0)$,  
then $(u_1,1)g\in F$ and $w g^{-1}(u_1,1)\in K(1)_{(u_1,1)}$. Thus
$w\in W^1 F\oplus W^2 F \oplus U^3 F$.  If $w_1=0$ then
 $w_2\neq 0$ and using a same argument as before we conclude that 
$(1,u_2)g\in F$, thus $wg^{-1}(1,u_2)\in K(1)_{(1,u_2)}$.  If $\widetilde{U}^3=
\widetilde{W}^1=0$ then $wg^{-1}(1,u_2)\in W^2$ and $w\in W^2 F$.
If some of the vector spaces $\widetilde{U}^3$, $\widetilde{W}^1$ are not null
then $u\in F$, from which we deduce that $g^{-1}(u_1,1)\in F$ and
 $wg^{-1}(u_1,1)\in  K(1)_{(u_1,1)}$. Hence $w\in W^1 F\oplus W^2 F \oplus U^3 F$.
\medbreak

If $S=\{b_i\}$ is any
basis of $V$ then $H$ is generated as an algebra by the set $$\{[b_i], g: b_i\in B, g\in G\}.$$ Indeed,
take $v\in V_1$, $w\in V_2$ then $(v,0)=\sum_i \alpha_i\, b_i$, $(0,w)=\sum_i \beta_i \, b_i$ for some
families of scalars $ \alpha_i,   \beta_i\in\ku$, then
$v=\sum_i \alpha_i\, [b_i]$ and  $w=\sum_i \beta_i \, [b_i] u$.
Let $\{b_i: i=1,\dots,r\}$ be a basis of $W=W^1\oplus W^2\oplus W^3$ and extend it to a basis $\{b_i: i=1,\dots,t\}$, $r\leq t$, of $V$.
Let $n> 1$ and $k\in K(n)$. Write 
$$k=\sum_{s_j\in \{0,1\},g_i\in G} \alpha_{s_1,\dots,s_t,i}\, [b_1]^{s_1} [b_2]^{s_2}\dots [b_t]^{s_t} g_i, $$
for some $ \alpha_{s_1,\dots,s_t,i}\in\ku$. Let $p:H\to H(1)$ be the canonical projection. Then
$(\id\ot p)\Delta(k)\in H(n-1)\otk K(1)$. It follows from a straightforward computation that
$ (\id\ot p)\Delta(k)$ is equal to
$$ \sum_l \sum_{s_j\in \{0,1\},g_i\in G} \alpha_{s_1,\dots,s_t,i} \big( \, h_{s_1,\dots,s_t,i}\ot [b_l] g_i +  \widetilde{h}_{s_1,\dots,s_t,i}\ot [b_l] u g_i\big),$$
for some $0\neq h_{s_1,\dots,s_t,i}, \widetilde{h}_{s_1,\dots,s_t,i} \in H(n-1)$. This implies that if $r< l$ and $s_l=1$ then 
$ \alpha_{s_1,\dots,s_t,i}=0$. Thus $K$ is generated as an algebra by $K(0)$ and $K(1)$,
whence $K=C(W^1,W^2,W^3, F)$.
\epf

\begin{defi}If $(W^1,W^2,W^3, F)$ is a coideal subalgebra data 
denote $\widetilde{C}(W^1,W^2,W^3, F)$ the subalgebra
of $\widetilde{H}$ generated by $\ku F$ and elements in
$W^1\oplus W^2$ and $\{[w]: w\in W^3\}$.
\end{defi}

\begin{teo}\label{coideal-tensor-prod-supergroup-tw} Let $K=\oplus_{i=0}^m K(i) \subseteq \widetilde{H}$
 be a homogeneous left coideal subalgebra. There
exists a coideal subalgebra data $(W^1,W^2,W^3, F)$ such that
$K=\widetilde{C}(W^1,W^2,W^3, F)$.
\end{teo}
\pf The proof follows the same argument as in the proof of Theorem 
\ref{coideal-tensor-prod-supergroup}.
\epf
\section{Module categories over tensor product of supergroup algebras}\label{section:modca}

We shall use the same notation as in the previous section, so we have a data $(V_1,V_2,u_1,u_2,G_1,G_2)$
 as in subsection \ref{tp-spa},  $H=\Ac(V_1,V_2,u_1,u_2,G_1,G_2)$
and $\widetilde{H}=\Hc(V_1,V_2,u_1,u_2,G_1,G_2)$. Denote $G=G_1\times G_2$, $H_i=\Ac(V_i,u_i,G_i)$ and 
$u=(u_1,u_2)\in G$.

We shall define  a family of comodule algebras over $H$ that will parameterize exact module
categories over $\Rep(H)$.

\begin{defi} We say that  the collection $(W, \beta, F, \psi)$
is a \emph{compatible data } with $(V_1,V_2,u_1,u_2,G_1,G_2)$ if
\begin{itemize}
 \item[(i)] $W=W^1\oplus W^2\oplus W^3 $ is a subspace of $V$ such that $(W^1,W^2,W^3, F)$ is a coideal subalgebra data,
 \item[(ii)] $\beta:W\times W\to \ku$
is a bilinear form stable under the action of $F$, such that
$$\beta(w_1,w_2)=-\beta(w_2,w_1),\; \beta(w_1,w_3)=\beta(w_3,w_1),\; \beta(w_2,w_3)=-\beta(w_3,w_2), 
$$
for all $w_i\in W^i$, $i=1,2,3$, and $\beta$ restricted to $W^i\times W^i$ is symmetric
for any $i=1,2,3$.
\item[(iii)] If $u\notin F$ then  $\beta$ restricted to $W^1\times W^2$ and $ W^2\times  W^3$ is null. 
\item[(iv)]  $\psi\in Z^2(F,\ku^{\times})$.
\end{itemize}

\end{defi}

Given a compatible data $(W, \beta, F, \psi)$ define
$\kc(W, \beta, F, \psi)$ as the algebra generated by $W$ and  $\{e_f: f\in F\}$,
subject to relations
$$ e_f e_h=\psi(f,h)\, e_{fh}, \quad e_f w= (f\cdot w)  e_f,$$
$$w_iw_j+w_iw_j=\beta(w_i,w_j) 1,\quad w_i\in W^i, w_j\in W^j,$$ 
for any $ (i,j)\in 
\{(1,1), (2,2), (1,3), (3,3) \},$ and relations
$$ w_2 w_3 - w_3 w_2= \beta(w_2,w_3)\, e_u , \text{ for any  } w_2\in W^2, w_3\in W^3,$$
$$ w_1 w_2 - w_2 w_1= \beta(w_1,w_2)\, e_u, \text{ for any  } w_1\in W^1, w_2\in W^2.$$
Define $\lambda: \kc(W, \beta, F, \psi)\to H\otk \kc(W, \beta, F, \psi)$ on the generators
$$\lambda(e_f)=f\ot e_f, \quad \lambda(w_1)= w_1\ot 1 + (u_1,1)\ot w_1, \text{
for all } f\in F, w_1\in W^1,$$
$$\lambda(w_2)= w_2\ot 1 + (1,u_2)\ot w_2 \quad \text{
for all } w_2\in W^2,$$ 
$$\lambda(v,w)=v\ot 1 + w (1,u_2)\ot e_u +  (u_1,1) \ot  (v,w), \quad \text{
for all } (v,w)\in W^3.$$

\begin{rmk}  If $(W, \beta, F, \psi)$ is a compatible data then $W$ comes with a
distinguished decomposition $W=W^1\oplus W^2\oplus W^3$. To be more precise
one should denote the algebras $ \kc(W, \beta, F, \psi)$ by 
$ \kc(W^1, W^2, W^3, \beta, F, \psi)$. We shall do this only in case  we want to emphasize
the direct decomposition of $W$.
\end{rmk}
\begin{defi}\label{def-ele} If $(0,0,W,F)$ is a coideal
subalgebra data and $(W, \beta, F, \psi)$
is a compatible data with $(V_1,V_2,u_1,u_2,G_1,G_2)$, we shall denote
$$\ele(W, \beta, F, \psi)= \kc(0,0,W, \beta, F, \psi).$$
\end{defi}

The algebras $\ele(W, \beta, F, \psi)$ will be the relevant ones when computing the Braur-Picard
group. 

\begin{prop}\label{right-simp} If $(W, \beta, F, \psi)$ is a compatible data then $\kc(W, \beta, F, \psi)$
 is a right $H$-simple left
$H$-comodule algebra with trivial coinvariants.
Also $\gr \Kc(W,\beta,F,\psi)=\Kc(W,0,F,\psi).$
\end{prop}
\pf The proof that these algebras are comodule algebras is straightforward. Also, it follows
from a direct computation that 
$$\gr \Kc(W,\beta,F,\psi)=\Kc(W,0,F,\psi),$$
and $\Kc(W,\beta,F,\psi)_0=\ku_\psi F$. Thus, the fact that these  algebras
are right $\Ac(V_1,V_2,u_1,u_2,G_1,G_2)$-simple follows from \cite[Prop. 4.4]{M1}.
\epf

Recall that in Section \ref{tib} we have defined  a left $H_i\otk H_i^{\cop}$-comodule
algebra $\diag(H_i)$. It follows from Lemma \ref{iso-cop} that there  is an isomorphism of Hopf algebras 
$$H_i\otk H_i^{\cop}\simeq H_i\otk H_i\simeq \Ac(V_i,V_i,u_i,u_i,G_i,G_i).$$
For any $i=1,2$ we shall denote $B_i=\Ac(V_i,V_i,u_i,u_i,G_i,G_i).$ Also
$$\diag(V_i)=\{(v,v)\in V_i\oplus V_i: v\in V_i\},$$
$$ \diag(G_i)=\{(g,g)\in G_i\times G_i: g\in G_i\}.$$

\begin{lema}\label{sigma-isom} For  any $i=1,2$ there is an isomorphism of left $B_i$-comodule algebras
$$\diag(H_i)\simeq \Kc(0\oplus 0 \oplus \diag(V_i),0,\diag(G_i),1).$$
\end{lema}
\pf Define $\sigma: \diag(H_i)\to \Kc(0\oplus 0 \oplus \diag(V_i) ,0,\diag(G_i),1)$ as follows.
For all $v\in V_i$, $g\in G_i$
$$\sigma(v)=(v,v)(u_i,u_i), \quad \sigma(g)=(g,g).$$
This gives a well-defined algebra isomorphism. It follows straightforward that
$\sigma$ is a $B_i$-comodule map.
\epf

\begin{rmk} In Lemma \ref{sigma-isom}  we write the space 
$W=0\oplus 0 \oplus \diag(V_i)$  to emphasize that $W^1=0, W^2=0$ and
$W^3= \diag(V_i)$.
\end{rmk}

\begin{prop}\label{comod=coid} Let $(W, 0, F, \psi)$ be a compatible data and $\widehat{\psi}\in Z^2(G,\ku^{\times})$ be a 2-cocycle such that
$\widehat{\psi}\mid_F=\psi$. Let $\sigma:H\otk H\to \ku$ be a Hopf 2-cocycle such that
 $\sigma(x,y)=\widehat{\psi}(x,y) $ for all $x,y\in G$, as defined in \eqref{lift-twist-def}.
Denote $\xi=\widehat{\psi}((u_1,1),(1,u_2))\widehat{\psi}((1,u_2),(u_1,1))^{-1}$.
If $\xi=1$ there is an isomorphism of comodule algebras
$$ \kc(W, 0, F, \psi) \simeq C(W^1,W^2,W^3, F)_\sigma.$$
If $\xi=-1$ there is an isomorphism of comodule algebras
$$ \kc(W, 0, F, \psi) \simeq \widetilde{C}(W^1,W^2,W^3, F)_\sigma.$$
\end{prop}
\pf One can  verify that the relations that hold in $C(W^1,W^2,W^3, F)_\sigma$
are the same relations in $ \kc(W, 0, F, \psi)$. Thus there is a well-defined projection
$C(W^1,W^2,W^3, F)_\sigma\twoheadrightarrow \kc(W, 0, F, \psi)$ which is an
isomorphism since both algebras have the same dimension.\epf 

The above Proposition can be extended when the bilinear form $\beta$ is not null. Let
us begin by constructing a Hopf 2-cocycle in $H$.
\begin{lema} Let $\sigma:H\otk H\to \ku$ be the 
map defined by
 $$\sigma(x,y)=\begin{cases} \widehat{\psi}(x,y) \text{ if } x,y\in G \\

               \frac{1}{2} \widehat{\psi}(sg,th)\beta(v,w)  \text{ if } x= vg, y=wh, v\in V_g, w\in
V_t, g,h\in G\\
0 \quad \text{ otherwise. }
               \end{cases}$$
Then $\sigma$ is a Hopf 2-cocycle.\qed
\end{lema}
Let $\xi=\widehat{\psi}((u_1,1),(1,u_2))\widehat{\psi}((1,u_2),(u_1,1))^{-1}$
and let $(W, \beta, F, \psi)$ be a compatible data.
\begin{prop}\label{twisting-comod-a}  If
$\xi=1$
then there is an isomorphism
of comodule algebras
$$ \kc(W, \beta, F, \psi) \simeq C(W^1,W^2,W^3, F)_\sigma .$$
If $\xi=-1$ there is an isomorphism of comodule algebras
$$ \kc(W, \beta, F, \psi) \simeq \widetilde{C}(W^1,W^2,W^3, F)_\sigma.$$
\end{prop}
\pf One can verify that the relations that hold in $C(W^1,W^2,W^3, F)_\sigma $
are the same relations that hold in $ \kc(W, \beta, F, \psi)$. Let us do this only for
$w_2\in W^2, (v,w)\in W^3$. By definition of $\sigma$ we have
$$w_2\cdot_\sigma [(v,w)]=\frac{1}{2}\beta(w_1,v) 1+\frac{1}{2} \beta(w_1,w) u +\widehat{\psi}(u_2,u_1) w_2 [(v,w)], $$
$$ [(v,w)] \cdot_\sigma w_2= \frac{1}{2}\beta(v,w_1) 1+\frac{1}{2} \beta(w,w_1) u+\widehat{\psi}(u_1,u_2)  [(v,w)]w_2. $$ 
Then 
$$w_2\cdot_\sigma [(v,w)] - [(v,w)] \cdot_\sigma w_2= \beta(w_1,w) u.$$
Thus there is a well-defined projection $C(W^1,W^2,W^3, F)_\sigma\twoheadrightarrow \kc(W, \beta, F, \psi)$ which is an
isomorphism since both algebras have the same dimension.
\epf

\begin{teo}\label{mod-over-supalg} Let $(V_1,V_2,u_1,u_2,G_1,G_2)$
be a data as in subsection \ref{tp-spa} and $H=\Ac(V_1,V_2,u_1,u_2,G_1,G_2)$. Let $\Mo$ be an indecomposable exact
left $\Rep(H)$-module category. Then there is a compatible data $(W,\beta,F,\psi)$ such that
 $\Mo$ is equivalent to the category ${}_{\Kc(W,\beta,F,\psi)}\Mo$ as $\Rep(H)$-modules.
\end{teo}
\pf
By Proposition \ref{right-simp} and \cite[Prop. 1.20]{AM} the families  ${}_{\Kc(W,\beta,F,\psi)}\Mo$ are exact indecomposable
module categories over $\Rep(H)$.
\medbreak

Let $\Mo$ be an indecomposable exact $\Rep(H)$-module category. Then, by \cite[Thm 3.3]{AM}
there exists a right $H$-simple left comodule algebra with trivial coinvariants $(A,\lambda)$
such that $\Mo={}_A\Mo$ as $\Rep(H)$-modules. Since $H$ is coradically graded then
$\gr A$ is a right $H$-simple left comodule algebra also with trivial coinvariants.

Since $H_0=\ku G$ and $A(0)$ is a left $\ku G$-comodule algebra right $\ku G$-simple
then there exists a subgroup $F\subseteq G$ and $\psi\in Z^2(F,\ku^{\times})$ such that
$A(0)=\ku_\psi F$.

Let $\widehat{\psi}\in Z^2(G,\ku^{\times})$ be a 2-cocycle such that
$\widehat{\psi}\mid_F=\psi$. Let $\sigma:H\otk H\to \ku$ be a Hopf 2-cocycle such that
 $\sigma(x,y)=\widehat{\psi}(x,y) $ for all $x,y\in G$, as defined in \eqref{lift-twist-def}.

By Lemma \ref{transp-comod} the algebra $(\gr A)_{\sigma^{-1}}$ is isomorphic to a homogeneous left 
coideal subalgebra of $H^{[\sigma^{-1}]}$. Set $\xi=\widehat{\psi}((u_1,1),(1,u_2))\widehat{\psi}((1,u_2),(u_1,1))^{-1}.$
Since $\xi^2=1$ then $\xi=\pm 1$. We shall analize what happends in both cases.

\subsection*{Case $\xi= 1$} It follows from Proposition  \ref{twisting-supergroup} that
 there is an isomorphism
of Hopf algebras $H^{[\sigma^{-1}]}\simeq \Ac(V_1,V_2,u_1,u_2,G_1,G_2)$, therefore 
$(\gr A)_{\sigma^{-1}}$ is isomorphic as a comodule algebra to 
a coideal subalgebra of $H$. Hence, from Theorem \ref{coideal-tensor-prod-supergroup}
we deduce that 
$(\gr A)_{\sigma^{-1}}=C(W^1,W^2,W^3, F)$ for some coideal subalgebra data
$(W^1,W^2,W^3, F)$. Proposition  \ref{comod=coid} 
implies that $(\gr A)\simeq \Kc(W,0,F,\psi)$.  Now, we have to determine all
liftings of $\Kc(W,0,F,\psi)$, that is all comodule algebras $A$ such that $(\gr A)\simeq \Kc(W,0,F,\psi)$.
\smallbreak

For any $w_1\in W^1$, $w_2\in W^2$, $(v,w)\in W^3$ let be $a_{w_1},a_{w_2},
 a_{(v,w)}\in A_1$  elements such that
$$\lambda(a_{w_1})=w_1\ot 1+ (u_1,1)\ot a_{w_1},\quad 
\lambda(a_{w_2})=w_2\ot 1+ (1,u_2)\ot a_{w_2},$$
$$\lambda(a_{(v,w)})= v\ot 1 + wu\ot e_u + (u_1,1)\ot a_{(v,w)},$$
 and the class of $a_w$ in $A(1)=A_1/A_0$ equals $w$. We can choose
these elements so that they satisfy that   
$$a_{v+w}= a_v+a_w, \quad f a_w f^{-1}= a_{f\cdot w}\quad \text{ for all } f\in F, v,w\in W.$$
The proof of the existence of such  elements is the same 
as the proof of \cite[Lemma 5.5]{M1}. Then $A$ is generated as an algebra by elements
$\{a_w , f: w\in W, f\in F\}$.

For any $ (i,j)\in \{(1,1), (2,2), (1,3),  (3,3) \}$ take $w_i\in W^i, w_j\in W^j$. 
Then
$$\lambda(a_{w_i}a_{w_j}+ a_{w_j}a_{w_i})= 1\ot a_{w_i}a_{w_j}+ a_{w_j}a_{w_i},$$
hence there exists an scalar $\beta(w_i, w_j)\in \ku$ such that
$$ a_{w_i}a_{w_j}+ a_{w_j}a_{w_i}= \beta(w_i, w_j)\, 1.$$
If $w_1\in W^1$, $w_2\in W^2$ then
$$\lambda(a_{w_1}a_{w_2}- a_{w_2}a_{w_1})= u\ot  a_{w_1}a_{w_2}- a_{w_2}a_{w_1},$$
hence there exists $\beta(w_1,w_2)\in \ku$ such that
$$a_{w_1}a_{w_2}- a_{w_2}a_{w_1} = \beta(w_1,w_2)\, e_u.$$
If $u\notin F$ then $ \beta(w_1,w_2)=0$.  The same is done in
the case $w_2\in W^2$, $w_3\in W^3$. One can prove that
$(W,\beta,F,\psi)$ is a compatible data and there is a comodule algebra projection
$$\kc(W,\beta,F,\psi) \twoheadrightarrow A$$
which is injective since both algebras have the same dimension.
\subsection*{Case $\xi=- 1$} The proof of this case is entirely similar to the case
 $\xi= 1$.
\epf

\subsection{Equivalence classes of module categories} We shall explain when two module
categories appearing in Theorem \ref{mod-over-supalg} are equivalent. 
\medbreak

Let $H$ be a finite-dimensional pointed Hopf algebra and $\Ac, \Ac'$ be right
$H$-simple left $A$-comodule algebras with trivial coinvariants. If $g\in H$ is a group-like
element we can define a new  comodule algebra $\Ac^g$ on
the same underlying algebra $\Ac$ with coaction given by $\lambda^g:\Ac^g \to H\otk \Ac^g,$
$\lambda^g(a)=g a\_{-1} g^{-1}\ot a\_0$, for all $a\in \Ac$.

\begin{teo}\label{Morita-equivalence}\cite[Thm. 4.2]{GM} The algebras $\Ac$, $\Ac'$
are equivariantly Morita equivalent  if and only if there exists an element
$g\in G(A)$ such that $\Ac'\simeq \Ac^{g}$ as comodule algebras.\qed
\end{teo}

\begin{teo}\label{Morita-equivalence2} Let $(V_1,V_2,u_1,u_2,G_1,G_2)$
be a data as in subsection \ref{tp-spa} and set $H=\Ac(V_1,V_2,u_1,u_2,G_1,G_2)$. 
Let $(W,\beta,F,\psi)$, $(U,\beta',F',\psi')$ be two compatible data. The module
categories  ${}_{\Kc(W,\beta,F,\psi)}\Mo$, ${}_{\Kc(U,\beta',F',\psi')}\Mo$
are equivalent if and only if there exists $g\in G$ such that
$$W^1=g\cdot U^1,\;   W^2=g\cdot U^2, \;W^3=g\cdot U^3,\;\beta'=g\cdot \beta, \;  F'=gFg^{-1},\; \;  \psi'=\psi^g.$$
Here $g\cdot \beta(v,w)= \beta(g^{-1}\cdot v,g^{-1}\cdot w)$ for all
$v,w\in U$.

\end{teo}

\pf
Let us prove that  if 
$\Kc(W,\beta,F,\psi)$ and $\Kc(W',\beta',F',\psi')$ are isomorphic as $H$-comodule algebras
then $W=W'$, $\beta=\beta'$, $F=F'$ and $\psi=\psi'$.
Let $\vartheta:\Kc(W,\beta,F,\psi)
\to \Kc(W',\beta',F',\psi')$ be an isomorphism of $H$-comodule algebras, then 
for any $f\in F$ we have that
$$ f\ot \vartheta(e_f)= \lambda( \vartheta(e_f)).$$
This implies that $ \vartheta(e_f)\in \Kc(W',\beta',F',\psi')_0=\ku F'$ and has no other 
possibility that being equal to $e_f$. Hence $F\subseteq F'$. The other inclusion can be
proven using the inverse of $\vartheta$.
Since $ \vartheta$ is an algebra morphism we deduce
that $\psi=\psi'$.

\medbreak

It is not difficult to prove that for any $i=1,2,3$ we have that $\vartheta(W^i)\subseteq U^i$. 
Since  $\vartheta$ is an isomorphism then $W^i=  U^i$ for any $i=1,2,3$. 
Since $ \vartheta$ is an algebra morphism the bilinear forms $\beta, \beta'$ must
be equal.

For any $g\in G$ there is an isomorphism of comodule algebras
$$\Kc(W,\beta,F,\psi)^g \simeq \Kc(g\cdot W,g\cdot \beta,gFg^{-1},\psi^g).$$
Indeed, the algebra map $\theta:\Kc(W,\beta,F,\psi)^g \to \Kc(g\cdot W,g\cdot \beta,gFg^{-1},\psi^g)$
determined by
$$\theta(w)=g\cdot w, \quad \theta(e_f)=e_{gfg^{-1}},$$
for all $w\in W$, $f\in F$, is an isomorphism of comodule algebras. The proof of the
Theorem follows now from Theorem \ref{Morita-equivalence}.
\epf
\section{The Brauer-Picard group of supergroup algebras}\label{sectio:brp}

The\emph{ Brauer-Picard groupoid} \cite{ENO} $\underline{\underline{\text{BrPic} }}$ 
is the 3-groupoid whose objects are finite tensor categories, 1-morphisms from
$\ca_1$ to $\ca_2$ are invertible exact $(\ca_1,\ca_2)$-bimodule categories, 2-morphisms are equivalences
 of such bimodule categories, and 3-morphisms are isomorphisms
of such equivalences. Forgetting the 3-morphisms and the 2-morphisms and identifying 1-morphisms one obtains
the groupoid BrPic. For a fixed tensor category $\ca$, the group $\text{BrPic}(\ca)$ consists of
equivalence classes of invertible exact $\ca$-bimodule categories and it is called the \emph{ Brauer-Picard group
of } $\ca$.

\medbreak

In this section $G$ will denote a finite Abelian group, $V$ is a finite-dimensional
$G$-module and $u\in G$ is an element of order 2
such that it acts on $V$ as $-1$. Also $H=\Ac(V,u,G)$. 

\subsection{The Brauer-Picard group of group algebras} Let us recall the results
obtained in \cite{ENO} on the computation of the Brauer-Picard group
of the category of representations of a finite Abelian group.

\medbreak

\begin{defi}\label{orth} Let $G$ be a finite Abelian group. The group $O(G\oplus \widehat{G})$ consists of group
 isomorphisms $\alpha:G\oplus \widehat{G}\to G\oplus \widehat{G} $ 
such that $\langle\alpha_2(g,\chi), \alpha_1(g,\chi) \rangle = \langle\chi, g\rangle$
for all $g\in G, \chi\in \widehat{G}$. Here $\alpha(g,\chi)= (\alpha_1(g,\chi) , \alpha_2(g,\chi))$.
 
\end{defi}

\begin{teo}\cite[Corollary 1.2]{ENO} There is an isomorphism of groups 
$$\text{BrPic}(\Rep(G))\simeq O(G\oplus \widehat{G}).$$\qed
\end{teo}

Let us explain how to obtain invertible bimodule categories
from elements in $O(G\oplus \widehat{G})$. Let $\alpha\in O(G\oplus \widehat{G})$ and define $U_\alpha\subseteq G\times G$ the subgroup
$$ U_\alpha=\{(\alpha_1(g,\chi) ,g): g\in G, \chi\in \widehat{G}\}.$$
and the 2-cocycle $\psi_\alpha:U_\alpha\times U_\alpha\to \ku^{\times}$ defined by
$$ \psi_\alpha((\alpha_1(g,\chi) ,g), (\alpha_1(h,\xi) ,h))=\langle\alpha_2(g,\chi)^{-1} , 
\alpha_1(h,\xi)\rangle \langle \chi,h\rangle.$$
It was proved in \cite{ENO} that the bimodule categories 
${}_{\ku_{\psi_\alpha}U_\alpha}\Mo$ are invertible and any invertible bimodule category is equivalent to one of this form.
 Note that
$U_{\id}=\diag(G), \psi_{\id}=1$.

\begin{exa}\label{exa-brC2} If $G=\Z_p$ for some prime $p\in \Na$ then $O(G\oplus \widehat{G})$ is isomorphic to 
the dihedral group $\Di_{2(p-1)}$. In particular if $p=2$ then 
$O(\Z_2\oplus \widehat{\Z_2})\simeq \Z_2$. The only non-trivial
element in $O(\Z_2\oplus \widehat{\Z_2})$ is $\gamma: \Z_2\oplus \widehat{\Z_2}
\to \Z_2\oplus \widehat{\Z_2}$ given by
\begin{equation}\label{br-c2}
 \gamma(u^i, \chi^j)=(u^j, \chi^i),
\end{equation}
for $i,j=0,1$. Here $u$ is the generator of $ \Z_2$ and $\chi$ is the generator of $ \widehat{\Z_2}$.

\end{exa}

\subsection{Families of invertible bimodule categories}

In this section we present families of invertible $\Rep(H)$-bimodule categories.

\begin{defi} We shall denote by $\ere(V,u,G)$ the set of   
collections $(W,\beta,\alpha)$, where
\begin{itemize}
 \item[(i)] $W\subseteq V\oplus V$ is a subspace such that  $ W\cap V\oplus 0=0=W\cap 0\oplus V$,
\item[(ii)] $\alpha\in O(G\oplus \widehat{G})$ is an isomorphism such that
$(u,u)\in U_\alpha$,

\item[(iii)] $W$ is invariant under the action of $U_\alpha$,
\item[(iv)]  $\beta:W\times W\to \ku$ is a symmetric bilinear form invariant under the
action of $U_\alpha$.
\end{itemize}
 
\end{defi}

If $(W,\beta,\alpha)$, $( \widetilde{W}, \widetilde{\beta},\widetilde{\alpha})$ are elements
in $\ere(V,u,G)$ we define
\begin{equation}\label{product-ere} (W,\beta,\alpha)\bullet ( \widetilde{W}, \widetilde{\beta},\widetilde{\alpha})=
(W\bullet  \widetilde{W}, \beta\bullet \widetilde{\beta}, \alpha\widetilde{\alpha}),
\end{equation}

where $W\bullet  \widetilde{W}$ is the subspace of $V\oplus V$
consisting of elements $(v_1,w_1)$ such that there exists a (necessarily unique) $v_2\in V$ such that
$(v_1,v_2)\in W$, $(v_2,w_1)\in \widetilde{W}$. The bilinear form $\beta\bullet \widetilde{\beta}$ 
is defined by
$$\beta\bullet \widetilde{\beta}((v_1,w_1), (v'_1,w'_1))= \beta((v_1,v_2), (v'_1,v'_2)) +
\widetilde{\beta}((v_2,w_1),(v'_2,w'_1)),  $$
where $v_2, v'_2\in V_2$ are the unique elements such that $(v_1,v_2), (v'_1,v'_2)\in W$ and
$ (v_2,w_1),$ $(v'_2,w'_1)\in \widetilde{W}$. The action of $U_{ \alpha\widetilde{\alpha}}$
on $W\bullet  \widetilde{W}$ is given as follows. If $g\in G, \chi\in \widehat{G}$, $(v,w)
\in W\bullet  \widetilde{W}$ then
$$(\alpha_1(\widetilde{\alpha}(g,\chi)), g)\cdot (v,w)=(\alpha_1(\widetilde{\alpha}(g,\chi))
\cdot v, g\cdot w). $$

\begin{lema} If $(W,\beta,\alpha)$, $( \widetilde{W}, \widetilde{\beta},\widetilde{\alpha})$
then  $(W,\beta,\alpha)\bullet ( \widetilde{W}, \widetilde{\beta},\widetilde{\alpha})
\in  \ere(V,u,G)$. 
\end{lema}
\pf We will only prove that the bilinear form $\beta\bullet \widetilde{\beta}$ is 
invariant under the action of $U_{ \alpha\widetilde{\alpha}}$. The other properties are
straightforward. Let $(v_1,w_1), (v'_1,w'_1)
\in W\bullet \widetilde{W}$ and $(f,g)\in U_{ \alpha\widetilde{\alpha}}$, then
$ \beta\bullet \widetilde{\beta}((f,g)\cdot (v_1,w_1), (f,g)\cdot (v'_1,w'_1))$ is equal to
\begin{align*}&= \beta\bullet \widetilde{\beta}( (f\cdot v_1,g \cdot w_1),  (f\cdot v'_1,g\cdot w'_1))\\
&= \beta((f\cdot v_1, x\cdot  v_2), (f\cdot v'_1,x\cdot  v'_2)) +
\widetilde{\beta}((x\cdot  v_2, g\cdot w_1),(x\cdot v'_2, g\cdot w'_1))\\
&= \beta((v_1,v_2), (v'_1,v'_2)) +
\widetilde{\beta}((v_2,w_1),(v'_2,w'_1))\\
&= \beta\bullet \widetilde{\beta}((v_1,w_1), (v'_1,w'_1)).
\end{align*}
In the above equalities the element $x\in G$ is the unique such that 
 $(f,x)\in U_{ \alpha}$ and  $(x,g)\in U_{\widetilde{\alpha}}$, and 
 $v_2, v'_2\in V_2$ are the unique elements such that $(v_1,v_2), (v'_1,v'_2)\in W$ and
$ (v_2,w_1),$ $(v'_2,w'_1)\in \widetilde{W}$.
\epf

\begin{defi}\label{defi-equi} We  say that
$(W,\beta,\alpha)\sim (\widetilde{W},\widetilde{\beta},\widetilde{\alpha})$
if there exists an element $g\in G\times G$ such that
$$\widetilde{W}=g\cdot W,\; \;  \widetilde{\beta}=g\cdot \beta, \alpha=\widetilde{\alpha}.$$
\end{defi}

If $(W,\beta,U_\alpha,\psi_\alpha)$ is a compatible family 
for some $\alpha\in O(G\oplus \widehat{G})$ we shall denote
$$\kc(W,\beta,\alpha)=\kc(W,\beta,U_\alpha,\psi_\alpha), \quad \ele(W,\beta,\alpha)
=\ele(W,\beta,U_\alpha,\psi_\alpha).$$

\begin{teo}\label{ex-inv-bimod} Let $(W,\beta, \alpha)\in\ere(V,u,G)$
such that there exists 
$(\widetilde{W},\widetilde{\beta}, \widetilde{\alpha})\in \ere(V,u,G)$
such that 
\begin{equation}\label{sp-invt1}(W,\beta, \alpha)\bullet(\widetilde{W},\widetilde{\beta}, \widetilde{\alpha})
\sim (\diag(V),0,\id), 
\end{equation}
\begin{equation}\label{sp-invt2}(\widetilde{W},\widetilde{\beta}, \widetilde{\alpha})\bullet
 (W,\beta, \alpha) \sim (\diag(V),0,\id).
\end{equation}

Then the $\Rep(H)$-bimodule category ${}_{\ele(W,\beta, \alpha)}
\Mo$ is invertible.
\end{teo}
\pf The proof is a (more complicated) version of the proof of the
fundamental theorem for Hopf modules \cite{Mo}. If
$L=\ele(W,\beta,\alpha), K=\ele(\widetilde{W},\widetilde{\beta}, \widetilde{\alpha})$ we shall prove that
the categories $ \Mo(H,H,K,\overline{L})$, ${}_{\ele(\diag(V),0,\id)}\Mo$ are equivalent as 
bimodule categories.
\medbreak

Let us fix some notation. If $(v_1, w_1) \in W\bullet \widetilde{W}$ 
then there exists a unique
$v_2\in V_2$ such that $(v_1, v_2)\in W$, $(v_2,w_1)\in \widetilde{W}$.
We shall denote
 $$\iota_1(v_1, w_1)=(v_1, v_2),\quad \iota_2(v_1, w_1)=(v_2,w_1).$$ Analogously if
$(v_1, w_1) \in   \widetilde{W}\bullet W$ there exists a unique $v_2\in V_1$
such that $(v_1,v_2)\in  \widetilde{W}$ and $(v_2,w_1) \in W$. We shall denote
$$ \widetilde{\iota}_1(v_1, w_1)= (v_1, v_2),\quad
 \widetilde{\iota}_2(v_1, w_1)=(v_2,w_1).$$

From \eqref{sp-invt1}, \eqref{sp-invt2} it follows that there are elements
$(a,b),(g,h)\in G\times G$ such that 
\begin{itemize}
 \item[(i)] $W\bullet \widetilde{W}=(g,h)\cdot \diag(V),$
\item[(ii)] $  \widetilde{W}\bullet W=(a,b)\cdot\diag(V).$
\end{itemize}

 Denote $S=  \ele(\diag(V),0,\id)$. Let  
$$\phi: S\to L\otk K, \quad \phib: S\to K\otk L,$$
 be the algebra morphisms determined as follows.
If $w\in V$, $f\in G$ then
$$\phi(w,w)= \iota_1(g\cdot w,h\cdot w)\ot 1 + e_{(u,u)} \ot  \iota_2(g\cdot w,h\cdot w),$$
$$ \phi(e_{(f,f)})=e_{(f, \widetilde{\alpha}(f,1))}\ot e_{(\widetilde{\alpha}(f,1),f)}.$$
 If $v\in V$,
 $x\in G$ then
$$\phib(v,v)=\widetilde{\iota}_1(a\cdot v,b\cdot v)\ot 1 + e_{(u,u)}\ot \widetilde{\iota}_2(a\cdot v,b\cdot v),$$
$$\phib(e_{(x,x)})= e_{(x,\alpha(x,1))}\ot  e_{(\alpha(x,1),x)}.$$
\begin{claim} The maps $\phi, \phib$ are well-defined.
\end{claim}
\pf[Proof of Claim] One should prove that for all $v,w\in V$, $f,g\in G\times G$
\begin{equation}\label{def-phi1} \phi(w,w) \phi(v,v)+ \phi(v,v)\phi(w,w) = 0,
\end{equation}
\begin{equation}\label{def-phi2} \phi(e_f) \phi(e_g)= \phi(e_{fg}),
\end{equation}
\begin{equation}\label{def-phi3} \phi(e_f)  \phi(w,w)=  \phi(f\cdot (w,w))\phi(e_f).
\end{equation}
The verification of these equalities is straightforward. The same equations hold
for $ \phib$.\epf

Let us recall the isomorphism of $\Ac(V,V,u,u,G,G)$-comodule
algebras $\sigma:\diag(H)\to S$ presented in the proof of Lemma \ref{sigma-isom}. We shall use the notation
$$ \phi(s)=\phi^1(s) \ot \phi^2(s), \;\, \phib(t)=\phib^1(t) \ot \phib^2(t),$$
omitting the summation symbol, for all $s,t\in S$.

\begin{claim} If $s,t\in S$ then
\begin{equation}\label{phi-comod1} \pi_2(\phib^1(t)_{-1})\Ss^{-1}(\pi_2(\phib^2(t)_{-1}))
\ot \phib^1(t)_0\ot \phib^2(t)_0 = 1 \ot  \phib(t)
\end{equation}
\begin{equation}\label{phi-comod2}  \pi_2( \phi^2(s)_{-1} )\Ss^{-1} (\pi_2( \phi^2(s)_{-1}))
\ot \phi^1(s)_0\ot \phi^2(s)_0= 1\ot  \phi(s).
\end{equation}
\end{claim}
The proof follows by verifying that both equalities hold for the generators of
the algebra $S$ and using that both maps $\phi, \phib$ are algebra morphisms.
\medbreak

If $M\in \Mo(H,H,K,\overline{L})$ define $\pi_M:M\to M^{\co H}$ the map
$$\pi_M(m)=\phib^1(\sigma(\Ss(m\_{-1})))\cdot  m\_0\cdot \phib^2(\sigma(\Ss(m\_{-1}))).$$
It follows from \eqref{phi-comod1} that the image of $\pi_M$ is indeed inside $M^{\co H}$.
The space $M^{\co H}$ has a left $S$-action given by
$$ s\cdot m= \phi^2(s)\cdot m \cdot \phi^1(s),$$
for all $s\in S_1$, $m\in M^{\co H}$. It follows from  \eqref{phi-comod2} that 
this action is well-defined, that is, if  $s\in S$, $m\in M^{\co H}$
then $ s\cdot m\in M^{\co H}$.

\medbreak

Let $\Gc: \Mo(H,H,K,\overline{L}) \to {}_{S}\Mo$ and $\Fc:{}_{S}\Mo
\to \Mo(H,H,K,\overline{L}) $ be the functors defined as follows. If
$M\in \Mo(H,H,K,\overline{L})$, $N\in {}_{S}\Mo$ then
$$\Fc(N)=(L\otk K) \ot_{S} N, \quad \Gc(M)=M^{\co H}.$$ 
The structure of right $S$-module on $L\otk K$ is given via $\phi$. Both functors
 are bimodule functors, see \cite[Prop. 3.7]{M3}. These functors are in fact
the same (up to some minor modifications)  functors described in
Section \ref{tib}. For any $M\in\Mo(H,H,K,\overline{L})$
 define
$$\alpha_M: M\to (L\otk K) \ot_{S}  M^{\co H}, \quad \beta_M: (L\otk K) \ot_{S}  M^{\co H}\to
M,$$
$$\alpha_M(m)=\phib(\sigma(m\_{-1})) \ot \pi_M(m\_0), \; \beta_M(l\ot k\ot m)=k\cdot m \cdot l,$$ 
for all $m\in M^{\co H}$, $l\in L$, $k\in K$.
\begin{claim} The maps $\alpha_M$, $\beta_M$ are inverse of each other.
\end{claim}
\pf[Proof of claim] Let $m\in M$ then $\beta_M\circ \alpha_M(m)$ is equal to
\begin{align*} &=\phib^1(\sigma(m\_{-1})) \cdot  \pi_M(m\_0) \cdot 
\phib^2(\sigma(m\_{-1}))\\
&=\phib^1(\sigma(m\_{-1})\sigma(\Ss(m\_0\_{-1})))\cdot  m\_0\_0 \cdot  
\phib^2(\sigma(m\_{-1})\sigma(\Ss(m\_0\_{-1})))\\
&= \epsilon(m\_{-1}) m\_0=m.
\end{align*}
Let $m\in M^{\co H}$, $l\in L$, $k\in K$. Then $ \alpha_M\circ\beta_M(l\ot k\ot m) $ is
equal to
\begin{align*} &= \alpha_M(k\cdot m \cdot l)\\
&= \phib\sigma\pi_2(k\_{-1}\Ss^{-1}(l\_{-1}) ) \ot \pi_M(k\_0\cdot m \cdot l\_0)\\
&=  \phib\sigma\pi_2(k\_{-1}\Ss^{-1}(l\_{-1}) ) \ot\\
&\otimes \phib^1\sigma\pi_2(k\_0\_{-1}\Ss^{-1}(l\_0\_{-1}))k\_0\_0\cdot  m \cdot l\_0\_0 \phib^2\sigma\pi_2(k\_0\_{-1}\Ss^{-1}(l\_0\_{-1}))
\end{align*}
 Now, to prove that  
$ \alpha_M\circ\beta_M(l\ot k\ot m) = (l\ot k\ot m)$ it is enough to prove that
\begin{align}\label{phi-r1} \phib\sigma\pi_2(k\_{-1})
\ot  \phib^1\sigma\pi_2(k\_0\_{-1}) k\_0\_0 \ot  \phib^2\sigma\pi_2(k\_0\_{-1}) 
= k\ot 1\ot 1\ot 1,
\end{align}
\begin{align}\label{phi-r2} \phib\sigma\pi_2(\Ss^{-1}(l\_{-1}))\ot 
\phib^1\sigma\pi_2(l\_0\_{-1})\ot l\_0\_0 \phib^2\sigma\pi_2(\Ss^{-1}(l\_0\_{-1}))
=1\ot l \ot 1\ot 1.
\end{align}
Since $ \phib\sigma\pi_2$ is an algebra map, equations \eqref{phi-r1},
\eqref{phi-r2} can be verified on the generators of the algebras $L$ and $K$.\epf
In conclusion we have that $\Fc \Gc =\Id$. Let us prove that $\Gc \Fc =\Id$. For any $N\in {}_{S}\Mo$ we have an inclusion
$$ N \hookrightarrow \Gc(\Fc(N)), \quad n\mapsto 1\ot 1\ot n,$$
for all $n\in N$. Let $\Psi: {}_{S}\Mo \to vect_\ku$ be the functor defined by
$\Psi(N)=\Gc(\Fc(N))/ N$ for all $N\in {}_{S}\Mo$. The functor $\Psi$ is a $\Rep(H)$-module
functor. Indeed, define $c_{X, N}:\Psi(X\otk N)\to X\otk  \Psi(N)$ by
$$ c_{X, N}(l\ot k\ot x\ot n)=\pi_1(l\_{-1})\cdot x\ot l\_0\ot k\ot n,$$
for all $X\in \Rep(H)$, $N\in {}_{1}\Mo$, $l\in L$, $k\in K$, $x\in X$, $n\in N$. 
It follows straightforward that $(\Psi, c)$ is a module functor, thus it is exact. 
The full subcategory $\No$ of objects such that $\Psi(N)=0$ is a submodule category of
${}_{S}\Mo$. Since $\Psi(S)= 0$ and the category ${}_{S}\Mo$ is indecomposable,
 then $\No={}_{S}\Mo$, which implies that $\Gc \Fc =\Id$.

\epf

\begin{teo}\label{p3pr} If $\alpha, \widetilde{\alpha}$ are elements in $ O(G\oplus \widehat{G})$ 
and $K=\ele(\diag(V),0,\alpha),$ $L=\ele(\diag(V),0, \widetilde{\alpha})$,  then there is an equivalence
of bimodule categories
$$\Mo(H,H,K,\overline{L}) \simeq 
{}_{\ele(\diag(V),0,\alpha\widetilde{\alpha})}\Mo.$$
\end{teo}
\pf The proof of Theorem \ref{ex-inv-bimod} applies \emph{mutatis mutandis} to this case.

\epf

\subsection{The Brauer-Picard group of supergroup algebras}

The comodule algebras $\ele(W,\beta, \alpha)$ will be the relevant ones
when computing the Brauer-Picard group for the representations categories
of supergroup algebras. In Theorem \ref{ex-inv-bimod} we have seen that
the bimodule categories  ${}_{\ele(W,\beta,\alpha)}\Mo$ are invertibles if
the compatible data  $(W,\beta, \alpha)$ is invertible in some sense. We must prove
now that these categories are the only invertible bimodule categories and we have to
describe the tensor product between them. In view of Theorem
\ref{tensor-bimod-hopf} we need first to  investigate the cotensor product of two 
such comodule algebras.

\medbreak

Let  $(W,\beta, \alpha)$,
$(\widetilde{W},\widetilde{\beta}, \widetilde{\alpha})$ be  compatible data  with 
$(V, u, G)$. We shall further assume that
$W$ and $\widetilde{W}$ have decompositions 
$$W=0\oplus 0 \oplus W^3, \quad \widetilde{W}=0\oplus 0 \oplus \widetilde{W}^3.$$

Let $L=\ele(W,\beta,\alpha),$ $K=\ele(\widetilde{W},\widetilde{\beta}, \widetilde{\alpha})$. If $(v_1, w_1) \in W\bullet \widetilde{W}$ 
then there exists a unique
$v_2\in V_2$ such that $(v_1, v_2)\in W$, $(v_2,w_1)\in \widetilde{W}$.
We shall denote
 $$\iota_1(v_1, w_1)=(v_1, v_2),\quad \iota_2(v_1, w_1)=(v_2,w_1).$$
 Analogously if
$(v_1, w_1) \in   \widetilde{W}\bullet W$ there exists a unique $v_2\in V$
such that $(v_1,v_2)\in  \widetilde{W}$ and $(v_2,w_1) \in W$. We shall denote
$$ \widetilde{\iota}_1(v_1, w_1)= (v_1, v_2), \quad
 \widetilde{\iota}_2(v_1, w_1)=(v_2,w_1).$$

Let $p_1, p_2:V\oplus V\to V$ the canonical projections, so $p_1(v,w)=v$ and
$p_2(v,w)=w$ for all $(v,w)\in V\oplus V$. Abusing
of the notation we shall also denote by $p_1, p_2:G\times G\to G$ the canonical projections.

\begin{lema}\label{basis-b} Let $\{(w_i^1, w_i^2)\}_{i=1}^t$ be a basis of $W\bullet \widetilde{W}\subseteq V\oplus
V$. There exists a basis $\{v_i\}_{i=1}^n$ of $W$ and a basis $\{w_i\}_{i=1}^m$ of $\widetilde{W}$
such that $t\leq n$, $t\leq m$ and $p_1(v_i)=w_i^1$, $p_2(w_i)=w_i^2$ for any $i=1,\dots, t$.
\end{lema}
\pf For any $i=1,\dots, t$ there exists $t_i\in V$ such that $(w_i^1, t_i)\in W$, 
$(t_i, w_i^2)\in\widetilde{W}$. The sets $\{(w_i^1, t_i)\}_{i=1}^t$,
$\{(t_i, w_i^2)\}_{i=1}^t$ are linearly independent, thus we can extend both set to a basis
in their corresponding spaces.\epf

1ga78
The right $H$-comodule structure described in
\eqref{comod-op2} will be denoted by $\lambda_r:L\to
L\otk H$ and if $x\in W$ then
\begin{equation}\label{right-comod1} \lambda_r(x)=e_{(u,u)}\ot p_2(x) + x\ot 1.
\end{equation}
The left $H$-comodule structure described in
\eqref{comod-op3} will be denoted by $\lambda_l:  K
\to H\otk  K$
and if $y\in \widetilde{W}$ then
\begin{equation}\label{left-comod1} \lambda_l(y)=p_2(y)\ot 1 + u \ot y.
\end{equation}

The proof of the next Lemma
can be done using an inductive argument. 
\begin{lema}\label{claims-b-coact} If $x=x_1\dots x_n e_f\in L$
and $y= y_1\dots y_m e_h\in K$ are elements such that $x_1,\dots ,x_n \in W$,
$y_1,\dots, y_m\in \widetilde{W} $
$f\in U_\alpha$, $h\in  U_{\widetilde{\alpha}}$ then
\begin{equation}\label{claims-b-coact1}  \lambda_r(x)=
\sum_{\substack{\epsilon_i,\delta_i\in \{0,1\}\\ \epsilon_i+\delta_i=1 }} \alpha^{\epsilon}_{\delta}\;\, x_1^{\epsilon_1}\dots x_n^{\epsilon_n} e^{\delta_1+\dots+\delta_n}_u
e_f \ot p_2(x_1)^{\delta_1} \dots p_2(x_n)^{\delta_n} p_2(f),
\end{equation}
\begin{equation}\label{claims-b-coact2} \lambda_l(y)=
\sum_{\substack{\epsilon_i,\delta_i\in \{0,1\}\\ \epsilon_i+\delta_i=1 }} \zeta^{\epsilon}_{\delta}\;\, p_1(y_1)^{\epsilon_1}\dots p_m(y_m)^{\epsilon_m}
u^{\epsilon_1+\dots+\epsilon_m} p_2(h) \ot y_1^{\delta_1}\dots y_m^{\delta_m} e_h,
\end{equation}
where all coefficients $ \alpha^{\epsilon}_{\delta}, \zeta^{\epsilon}_{\delta}\in \ku$ are not null. \qed
\end{lema}

\begin{prop}\label{cotensor-k3}  
Assume that $\widetilde{\alpha}=\id$. Then there is an isomorphism of left $H\otk H^{\cop}$-comodule algebras
\begin{equation}\label{cotensor-ele} \ele(W,\beta, \alpha)\Box_{H} \ele(\widetilde{W},\widetilde{\beta}, \id)
\simeq \ele(W\bullet \widetilde{W}, \beta\bullet \widetilde{\beta},\alpha).
\end{equation}
\end{prop}
\pf    Note that $ \ele(\widetilde{W},\widetilde{\beta}, \id)_0=\diag(G)$. Let
$\{x_1,\dots, x_n\}$ be a basis of $W$ and
$\{y_1,\dots, y_m\}$ be a basis of $\widetilde{W}$ such that they are
extensions of a basis of $W\bullet\widetilde{W}$ in the sense of Lemma \ref{basis-b}. Without loss of generality we can assume
that $m\leq n$.
\medbreak

For any $0\leq s \leq n,
0\leq t\leq m$ define $L(s)$  the subspace of $L$ generated by elements of the form
$$x_1^{\epsilon_1}\dots x_n^{\epsilon_n} e_f, \quad \text{ where } \epsilon_i=0,1,\;\;  
\epsilon_1+\dots + \epsilon_n=s,\;\; f\in U_\alpha.$$
 Analogously, define $K(t)$ 
 the subspace of $K$ gernerated by elements of the form
$$y_1^{\delta_1}\dots y_m^{\delta_m} e_f, \quad \text{ where }
\delta_j=0,1,\;\; \delta_1 +\dots +\delta_m=t, \;\;f\in \diag(G).$$
Then  $L=\oplus_{s=0}^n L(s)$ and $K=\oplus_{t=0}^m K(t)$.

\medbreak

 For any $i=1,\dots, m$ set $w_i=(p_1(x_i), p_2(y_i))$,
$\pi:K \to \ku 1$ the canonical projection and define $S\subseteq F$ the subset
of elements $(f_1,f_2)$ such that there exists $(f_2,g)\in \widetilde{F}$. We shall denote by
$p:\ku F\to  \ku F\bullet \widetilde{F}$ the linear map determined by
$$p(e_{(f_1,f_2)})=\begin{cases} 0\quad \quad \quad \text{ if } (f_1,f_2)\notin S \\
       e_{(f_1,g)}  \quad \text{ if } (f_2,g)\in \diag(G).
       \end{cases}
$$
By the assumptions on $F$ and $\widetilde{F}$ the map $p$ is well-defined. 
Define the
map $\theta:L \Box_{H}K \to \ele (W\bullet \widetilde{W}, \beta\bullet \widetilde{\beta}, \alpha) $ as follows. If $d\in\Na$ and $z\in L \Box_{H}K$ is an element of the form
\begin{equation}\label{element-z} z= \sum_{\substack{a_1+\dots +a_m+ b_1+\dots +b_n=d\\ f\in U_\alpha,\, h\in  \diag(G)}}
 \beta^{a, f}_{b, h}\;\; x_1^{a_1}\dots x_n^{a_n} e_f\ot
y_1^{b_1}\dots y_m^{b_m} e_h,
\end{equation}
then
$$\theta(z)=\sum_{\substack{a_1+\dots +a_m=d\\ \\f\in U_\alpha}}
 \beta^{a, f}_{0, 1}\;\;  w_1^{a_1}\dots w_m^{a_m} e_{(p_1(f),p_2(h))}. $$
From now on we shall write $a+ b=d$ when we mean that $a_1+\dots +a_m+ b_1+\dots +b_n=d$.
\begin{claim}\label{theta-iso}  The map $\theta$ is a well-defined injective linear map. In particular we have that
$\dim(L \Box_{H}K)\leq \dim(\ele (W\bullet \widetilde{W}, \beta\bullet \widetilde{\beta}, \alpha)). $
\end{claim}

\pf[Proof of claim]
We must prove that  $\theta$ is well-defined in $L \Box_{H}K$ and that it is injective. Let us prove the first. Let be $0 \leq d\leq n+m$ and $z\in 
L \Box_{H}K \cap \oplus_{s=0}^{n+m}\; L(s)\otk K(d-s)$ a non-zero element, then there are scalars $ \beta^{a, f}_{b, h}\in  \ku$ such that
\begin{equation}\label{element-in-cot} z= \sum_{\substack{a+ b=d\\ f\in U_\alpha, h\in \diag(G)}}
 \beta^{a, f}_{b, h}\;\; x_1^{a_1}\dots x_n^{a_n} e_f\ot
y_1^{b_1}\dots y_m^{b_m} e_h.
\end{equation}
Using \eqref{claims-b-coact1}, \eqref{claims-b-coact2} one gets that
\begin{equation}\label{dist-cotensor1} \sum_{\substack{a+ b=d\\ \epsilon_i+ \delta_i=a_i\\ f\in U_\alpha, h\in  \diag(G)}}
 \beta^{a, f}_{b, h}  \alpha^{\epsilon}_{\delta} \;
 x_1^{\epsilon_1}\dots x_n^{\epsilon_n} e^{\epsilon}_u
e_f \ot p_2(x_1)^{\delta_1} \dots p_2(x_n)^{\delta_n} p_2(f) \ot
y_1^{b_1}\dots y_m^{b_m} e_h
\end{equation}
equals to
\begin{equation}\label{dist-cotensor2} \sum_{\substack{a+ b=d\\ \epsilon_i+ \delta_i=b_i\\ f\in U_\alpha, h\in  \diag(G)}} 
 \beta^{a, f}_{b, h} \zeta^{\epsilon}_{\delta}\;  x_1^{a_1}\dots x_n^{a_n} e_f\ot
 p_1(y_1)^{\epsilon_1}\dots p_m(y_m)^{\epsilon_m}
u_2^{\epsilon} p_2(h) \ot y_1^{\delta_1}\dots y_m^{\delta_m} e_h.
\end{equation}
Since $z\neq 0$ there exists some $ \beta^{a, f}_{b, h}\neq 0$. Define
$$I(z)=\{1 \leq i \leq n: \text{ there exists }  \beta^{a', f}_{b', h}\neq 0 \text{ and } a'_i=1\}.$$
Let us assume that $1\in I(z)$, thus there exists some $\beta^{a', f}_{b', h}\neq 0$
where $a'_1=1$. The next argument does not depend on this choice but it simplifies
the notation.

Comparing  elements \eqref{dist-cotensor1} and \eqref{dist-cotensor2} we conclude, perhaps after reordering the elements of
the basis $\{y_1,\dots,y_m\}$, that
\begin{equation} \sum_{\substack{a+ b=d\\f\in U_\alpha, h\in  \diag(G)}}
  \beta^{a, f}_{b, h}  \alpha^{(0,a_2,\dots,a_n)}_{(1,0,\dots,0)} \;
 x_2^{a_2}\dots x_n^{a_n} e^{a}_u
e_f \ot p_2(x_1)p_2(f) \ot
y_1^{b_1}\dots y_m^{b_m} e_h
\end{equation}
must be equal  to
\begin{equation} \sum_{\substack{a+ b=d\\ f\in U_\alpha, h\in  \diag(G)}} 
 \beta^{a, f}_{b, h} \zeta^{(0,a_2,\dots,a_n)}_{(1,0,\dots,0)}\;  x_1^{a_1}\dots x_n^{a_n} e_f\ot
 p_1(y_1)
u_2 p_2(h) \ot y_2^{b_2}\dots y_m^{b_m} e_h.
\end{equation}
Since $\beta^{a', f}_{b', h}\neq 0$ then $p_1(y_1)=p_2(x_1)$ and $ p_2(f)=u p_2(h)=p_2((u,u) h)$.
 
\medbreak

Let $\lambda_1: L\otk K\to H\otk H^{\cop}\otk  L\otk K$ be the coaction
given in \eqref{coact-coprodu-t}, that is if $l\ot k\in  L\otk K$ then
$$ \lambda(l\ot k)=\pi_{1}(l\_{-1})\ot \pi_{1}(k\_{-1}) \ot l\_0\ot k\_0.$$
With this coaction $L\otk K$ is a comodule algebra and has $L\Box_{H}  K$
is a subcomodule algebra. Taking $H'=H\otk H^{\cop}$, 
$A= L\otk K$ and $B=L\Box_{H}  K$ we are under the hypothesis of
Lemma \ref{comod-st-tech}. This implies that any element  $z\in L\Box_{H}  K$
can be written as
$$ z=\sum_{d=0}^{n+m}\; z_d$$
where $z_d\in L \Box_{H}K \cap \oplus_{s=0}^{n+m}\; L(s)\otk K(d-s).$ 
Let us prove now that $\theta$ is injective. Assume that $z$ is an element as
in \eqref{element-z} such that $\theta(z)=0$ and $z\neq 0$. Hence any 
$ \beta^{a, f}_{0, 1}=0$ for any $a$ such that $a_1+\dots +a_m=d$.
 Since $z\neq 0$ there exists at least one coefficient $\beta^{a', f}_{b', h}\neq 0$. Let
us compare the coefficient of the term
\begin{equation}\label{comparing-coeff}
  x_1^{a'_1}\dots x_n^{a'_n} e_f\ot  p_1(y_1)^{b'_1}\dots p_m(y_m)^{b'_m}
u^{b'_1+\dots b'_m}\ot 1
\end{equation}
in equations \eqref{dist-cotensor1},  \eqref{dist-cotensor2}. The coefficient of the term \eqref{comparing-coeff}
in the summand \eqref{dist-cotensor1} is $\beta^{\widetilde{a}, f}_{0, 1} \alpha^{a'}_{b'}$, for
some $\widetilde{a}=(\widetilde{a}_1,\dots, \widetilde{a}_n)$ such that $\widetilde{a}_1+\dots
+\widetilde{a}_n=d$,
and  in the summand \eqref{dist-cotensor2} is $\beta^{a', f}_{b', 1} \zeta^{b'}_{0}$.
Thus $\beta^{a', f}_{b', 1} \zeta^{b'}_{0}= \beta^{\widetilde{a}, f}_{0, 1} \alpha^{a'}_{b'}=0$,
whence $\beta^{a', f}_{b', h}=0$, which is a contradiction, thus $\theta$ is injective. This finishes the
proof of the claim. 
\epf
Define
$\phi:\ele (W\bullet \widetilde{W}, \beta\bullet \widetilde{\beta}, \alpha)\to L\Box_{H} K$ the algebra map determined as follows.
If $w\in W\bullet \widetilde{W}$ then
$$\phi(w)= \iota_1(w)\ot 1 + e_u \ot  \iota_2(w),$$
and if $(f,g)\in  U_\alpha$ then 
$$\phi(e_{(f,g)})=e_{(f,g)}\ot e_{(g,g)}.$$

 The map $\phi$ extends to a comodule  algebra morphism and the image
is contained in $L\Box_{H} K$. To prove that $\phi$ is well-defined one should verify that
\begin{equation}\label{def-phi1} \phi(w) \phi(v)+ \phi(v)\phi(w) = \beta\bullet \widetilde{\beta}(v,w) 1,
\end{equation}
\begin{equation}\label{def-phi2} \phi(e_f) \phi(e_g)= \psi_{\alpha}(f,g)\; \phi(e_{fg})
\end{equation}
\begin{equation}\label{def-phi3} \phi(e_f)  \phi(w)=  \phi(f\cdot w)\phi(e_f),
\end{equation}
for all $w,v\in W\bullet \widetilde{W}$, $f,g\in U_\alpha$. This is done by a straightforward
computation. To prove that the image of $\phi$ is contained in $L\Box_{H} K$
we must prove that if
 $w\in W\bullet \widetilde{W}$ then $$\iota_1(w)\ot 1 + e_u \ot  \iota_2(w)\in 
L\Box_{H_2} K.$$ 
This calculation is readily proven. Let us prove that $\phi$ is a comodule morphism. For the moment we shall
denote by $\lambda_{W\bullet \widetilde{W}}$ the coaction of $\ele (W\bullet \widetilde{W}, \beta\bullet \widetilde{\beta}, \alpha)$. Let $w=(v_1,w_1)\in 
W\bullet \widetilde{W}$ then
\begin{align*}(\id\ot \phi)\lambda_{W\bullet \widetilde{W}}(w)&=
(\id\ot \phi)(v_1\ot 1+ w_1 (u_1,u_1)\ot e_{ (u_1,u_1)} + (u_1,1)\ot w)\\
&=v_1\ot 1\ot 1+ w_1 (u_1,u_1)\ot e_{ (u_1,u_2)} \ot e_{ (u_2,u_1)}+\\
& + (u_1,1)\ot \iota_1(w)\ot 1 + (u_1,1)\ot e_{ (u_1,u_2)} \ot  \iota_2(w).
\end{align*}
Let $\lambda$ denotes the coaction of 
$L \Box_{H}K$
described in \eqref{coact-coprodu-t} and 
$$\iota_1(w)=(v_1, v_2),\quad  \iota_2(w)=(v_2,w_1).$$
We have that
\begin{align*}\lambda (\iota_1(w)\ot 1)&=(\pi_{H_1}\ot \pi_{H_1}\ot \id)(v_1\ot 1\ot 1  + v_2 u\ot e_u\ot 1 
+ u_1\ot \iota_1(w)\ot 1 )\\
&= v_1\ot 1\ot 1 +  u_1\ot \iota_1(w)\ot 1,
\end{align*}
and
\begin{align*}\lambda (e_u \ot  \iota_2(w))&=(\pi_{H_1}\ot \pi_{H_1}\ot \id)\big(
uv_2\ot e_u \ot 1+u w_1 (u_2,u_1)\ot e_u\ot e_{(u_2,u_1)}+\\
& + u(u_2,1)\ot e_u\ot \iota_2(w)\big)\\
&=(u_1,1) w_1 (1,u_1)\ot e_u\ot e_{(u_2,u_1)}+(u_1,1) \ot e_u\ot \iota_2(w)\\
&= w_1 (u_1,u_1)\ot e_u\ot e_{(u_2,u_1)}+(u_1,1) \ot e_u\ot \iota_2(w)
\end{align*}
The last equality follows because $u_1$ commutes with $w_1$. Then 
$$ (\id\ot \phi)\lambda_{W\bullet \widetilde{W}}(w) = \lambda \phi(w).$$
An easy computation shows that the same equality holds for the group elements
in $U_\alpha$.  Clearly the map $\phi$ is injective. This implies that 
$ \dim(\ele (W\bullet \widetilde{W}, \beta\bullet \widetilde{\beta}, F\bullet \widetilde{F}
, \psi\bullet \widetilde{\psi}))\leq \dim(L \Box_{H_2}K),$  but from Claim
 \ref{theta-iso} it follows that both spaces have the same dimension. Therefore 
$\phi$ is an isomorphism.

\epf

 Let  $(W,\beta, F,\psi)$, $(\widetilde{W},\widetilde{\beta}, \widetilde{F},\widetilde{\psi})$ 
 be  compatible data  with 
$(V,V,u,u,G,G)$. The spaces $W, \widetilde{W}$
have  decompositions $W=W^1\oplus W^2 \oplus W^3$, 
$\widetilde{W}=\widetilde{W}^1 \oplus \widetilde{W}^2 \oplus\widetilde{W}^3. $

\medbreak 

Let 
$L=\kc(W,\beta, F,\psi),$ $K=\kc(\widetilde{W},\widetilde{\beta}, \widetilde{F},\widetilde{\psi})$.
The tensor product $L\otk K$ has a left $H$-comodule structure 
$\delta:L\otk K\to H_2\otk L\otk K$ given by
$$ \delta(l\ot k)=\pi_2(k\_{-1})\Ss^{-1}(\pi_2(l\_{-1}))\ot l\_0\ot k\_0,$$
for all $l\ot k\in L\otk K$. This coaction was already used in \eqref{comod-str-B}.
\begin{prop}\label{freeness} The following assertions hold. 
\begin{itemize}

 \item[1.]  If ${}_K\Mo$ is an invertible bimodule category then $\widetilde{W}^2=0$.
 \item[2.]  If ${}_K\Mo$ is an invertible bimodule category then $\widetilde{W}^1=0$.
 \item[3.]   If ${}_K\Mo$ is an invertible bimodule category then 
 $$\widetilde{F}=U_\alpha, \quad \widetilde{\psi}= \psi_\alpha,$$
for some  $\alpha\in O(G\oplus \widehat{G})$.

 \item[4.] If $W^2=W^1=\widetilde{W}^2=\widetilde{W}^1=0$
and $\widetilde{F}=U_{\widetilde{\alpha}},$ $\widetilde{\psi}= \psi_{\widetilde{\alpha}},$
$F=U_\alpha,  \widetilde{\psi}= \psi_\alpha$ for some
$\alpha, \widetilde{\alpha}\in O(G\oplus \widehat{G})$ there is an isomorphism  
$$L\otk K\simeq N\otk (L\Box_{H} K)$$ of right $L\Box_{H} K$-modules
and left $H$-comodules, where $N$ is a certain
left $H$-comodule.
\end{itemize}

\end{prop}
 \pf 1. Since 
${}_K\Mo$ is an invertible bimodule category then
$$\big( {}_K\Mo\big)^{\op}\boxtimes_{\Rep(H)} {}_K\Mo \simeq 
\Mo(H,H,K,\overline{K})\simeq \Rep(H_2).$$
For any vector space $X$ and $P\in \Mo(H,H,K,\overline{K})$ we write
$X\otk P$ the object in the category $\Mo(H,H,K,\overline{K})$ with structure
concentrated in $P$. Let $\No$ be the full  subcategory of $\Mo(H,H,K,\overline{K})$
consisting of objects $P$ such that $X\otb P\simeq X\otk P$ for all $X\in \Rep(H\ot H^{\cop})$.
The category $\No$ is a submodule category of $\Mo(H,H,K,\overline{K})$. 
It could not happend that $\No$ equals $\Mo(H,H,K,\overline{K})$ since
$\Mo(H,H,K,\overline{K})$ is
equivalent to $\Rep(H)$. Thus $\No$ must be the null category.

Let us assume that $\widetilde{W}^2\neq 0$. Let $<W^2>$ be the subalgebra
of $K$ generated by elements in $W^2$. We have inclusions 
$$S=\overline{<W^2>}\Box_{H}<W^2> \hookrightarrow \overline{K}\Box_{H_2}K
\hookrightarrow  \overline{K}\otk K,
$$
of left $H\ot H^{\cop}$-comodule algebras. Note that the coaction
of $S$ is trivial, that is, if $\delta:S\to H\ot H^{\cop}\otk S$,
$\sum k\ot l\in S$ then $\delta(\sum k\ot l)=1\ot 1\ot \sum k\ot l$. This implies that
for any $X\in \Rep(H\ot H^{\cop})$ and $M\in {}_S\Mo$
$X\otb M=X\otk M$, where the $S$-action  on $X\otk M$ is concentrated in the second tensorand.
From this observation we deduce that for any $M\in {}_S\Mo$ the object
$\overline{K}\otk K\ot_S M$ belongs to $\No$. This is a contradiction, which means that
$\widetilde{W}^2=0$.
\medbreak

2. It follows by using the same argument as in item (1).
\medbreak

3. Let us assume that ${}_L\Mo$ is the inverse of the bimodule
category ${}_K\Mo$. From the previous results we know that
$W^2=W^1=\widetilde{W}^2=\widetilde{W}^1=0$. Let us prove that
$(L\Box_H K)_0= L_0 \Box_{H_0} K_0$. The inclusion
$(L\Box_H K)_0\supseteq L_0 \Box_{H_0} K_0$ is immediate. Let $\sum l\ot k\in (L\Box_H K)_0$,
then
$$ \pi_1(l\_{-1})\ot l\_0 \in H_0\otk L, \quad  \pi_1(k\_{-1})\ot k\_0 \in H_0\otk K.$$
The only possibility for this to happen is that $l\in L_0, k\in K_0$. Now
the result follows from \cite[Corollary 5.6]{M3}.
\medbreak
4.   It follows from Proposition \ref{twisting-comod-a} and from \eqref{cotensor-ele}
 that  $L \Box_{H} K$ is a twisting $C_\sigma$ of some coideal subalgebra $C$
of either $\Ac(V,V,u,u,G,G)$ or $\Hc(V,V,u,u,G,G)$.
This means that there are equivalences of categories
$${}^B\Mo_{L \Box_{H} K}\simeq  {}^{B^{[\sigma]}}\Mo_{C}  \simeq{}^Q\Mo,$$
where $B$ is either $\Ac(V,V,u,u,G,G)$ or $\Hc(V,V,u,u,G,G)$
and $Q=B/BC^+$. The first equivalence is \cite[Lemma 2.1]{M1} and the second one is standard, see e.g. \cite{S}.
Thus any object of ${}^B\Mo_{L \Box_{H} K}$ is equivalent to
$B^{[\sigma]}\Box_Q N$ for some $N\in {}^Q\Mo$. In particular, since
$L\otk K\in {}^B\Mo_{L \Box_{H} K}$ there exists $N\in {}^Q\Mo$
such that $L\otk K\simeq B^{[\sigma]}\Box_Q N$. Since $B^{[\sigma]}\simeq C\otk Q$
as right $C$-modules and left $Q$-comodules, then  $L\otk K\simeq C_\sigma  \otk N$.
\epf

\begin{teo}\label{main-on-prod} If $\alpha, \widetilde{\alpha}$ are elements in $ O(G\oplus \widehat{G})$,  then there is an equivalence
of bimodule categories
$$ {}_{\ele(W,\beta,\alpha)}\Mo \boxtimes_{\Rep(H)}
 {}_{\ele(\widetilde{W},\widetilde{\beta}, \widetilde{\alpha})}\Mo \simeq 
{}_{\ele(W\bullet \widetilde{W},\beta\bullet \widetilde{\beta},\alpha\widetilde{\alpha})}\Mo.$$
\end{teo}
\pf From Proposition \ref{freeness} (4) we can apply  Theorem \ref{tensor-bimod-hopf}
and we get that
$${}_{\ele(W,\beta,\alpha)}\Mo\simeq {}_{\ele(W,\beta,\id)\Box_H \ele(V,0,\alpha)}\Mo
\simeq  {}_{\ele(W,\beta,\id)}\Mo  \boxtimes_{\Rep(H)}  {}_{\ele(V,0,\alpha)}\Mo, $$
where the first isomorphism is \eqref{cotensor-ele}. Then 
$${}_{\ele(W,\beta,\alpha)}\Mo \boxtimes_{\Rep(H)}
 {}_{\ele(\widetilde{W},\widetilde{\beta}, \widetilde{\alpha})}\Mo$$
is isomorphic to
$$  {}_{\ele(W,\beta,\id)}\Mo  \boxtimes_{\Rep(H)}  {}_{\ele(V,0,\alpha)}\Mo  \boxtimes_{\Rep(H)}  
{}_{\ele(V,0,\widetilde{\alpha})}\Mo \boxtimes_{\Rep(H)} {}_{\ele(\widetilde{W},\widetilde{\beta}, \id)}\Mo.$$
Using Theorem \ref{p3pr} we obtain that this tensor product is isomorphic to
$$  {}_{\ele(W,\beta,\id)}\Mo  \boxtimes_{\Rep(H)}  {}_{\ele(V,0,\alpha \widetilde{\alpha})}\Mo  \boxtimes_{\Rep(H)}  
 {}_{\ele(\widetilde{W},\widetilde{\beta}, \id)}\Mo,$$
and using again Theorem  \ref{tensor-bimod-hopf} we get the result.
\epf

Define $\nic( V, u,G)$ to be the group of invertible elements in $\ere(V,u,G)/ \sim$
with product $\bullet$ described in \eqref{product-ere}. 

\begin{teo} Let $G$ be a finite group, $u\in G$ be a central element of order 2 and $V$ a finite-dimensional 
$G$-module such that $u\cdot v=-v$  for all $v\in V$. There is an isomorphism of groups 
$$\brp(\Rep(\Ac(V,u,G)))\simeq \nic( V, u,G).$$
\end{teo}
\pf It follows from Theorem \ref{ex-inv-bimod} that the application 
$$\ere(V, u,G)^{\times}\to\brp(\Rep(\Ac(V,u,G))), \;  (W,\beta, \alpha)
\mapsto {}_{\ele(W,\beta, \alpha)}
\Mo$$
is well-defined. It follows from Theorem \ref{Morita-equivalence2} that  
$(W,\beta, \alpha)\sim (\widetilde{W},\widetilde{\beta}, \widetilde{\alpha})$
if and only if the module categories ${}_{\ele(W,\beta, \alpha)}
\Mo$, ${}_{\ele(\widetilde{W},\widetilde{\beta}, \widetilde{\alpha})}\Mo$ are equi\-valent. Hence we have a well-defined
injective map
$$\nic( V, u,G) \to\brp(\Rep(\Ac(V,u,G)), \;  (W,\beta, \alpha)\mapsto {}_{\ele(W,\beta, \alpha)}
\Mo.$$
Proposition \ref{cotensor-k3} implies that this map is a group homomorphism. Let
us prove that it is surjective. Let $\Mo$ be an exact invertible $\Rep(\Ac(V,u,G))$-bimodule
category. Then, by Theorem \ref{mod-over-supalg} there exists a data $(W^1\oplus W^2 \oplus W^3,\beta, F,\psi)$ compatible with $(V,V,u,u,G,G)$
and an equivalence $\Mo\simeq {}_{\kc(W^1,W^2,W^3,\beta, F,\psi)}
\Mo$ of bimodule categories. By Proposition \ref{freeness} (1) and (2) $W^1=W^2=0$. Also
by Proposition \ref{freeness} (3) there exists $\alpha\in O(G\oplus \widehat{G})$ such
that $(F,\psi)=(U_\alpha,\psi_\alpha)$, thus $\kc(W^1,W^2,W^3,\beta, F,\psi)=
\ele(W,\beta, \alpha)$.

Since $\Mo$ is invertible there exists another
compatible data $(\widetilde{W},\widetilde{\beta}, \widetilde{\alpha})$ such that 
\begin{equation}\label{tp111}
 {}_{\ele(W,\beta, \alpha)}
\Mo \boxtimes_{\Rep(\Ac(V,u,G))} {}_{\ele(\widetilde{W},\widetilde{\beta}, \widetilde{\alpha})}\Mo
\end{equation}
is equivalent to ${}_{\ele(\diag(V),0,\diag(G),1)}\Mo $. From Theorem \ref{main-on-prod}
we conclude  that the tensor product category \eqref{tp111} is equivalent to the 
category
$$ {}_{\kc( W\bullet \widetilde{W},
\beta\bullet \widetilde{\beta},\alpha\widetilde{\alpha} )}\Mo.$$
It follows from Theorem \ref{Morita-equivalence2} that
$( W\bullet \widetilde{W},
\beta\bullet \widetilde{\beta},\alpha\widetilde{\alpha} )\sim (\diag(V),0,\diag(G),1)$
and therefore $(W,\beta, \alpha)\in \nic( V, u,G)$. This finishes the proof of the Theorem.
\epf

\subsection{Another description of the Brauer-Picard group}
In  \cite{ENO} the authors give a beautiful description of the group
$\text{BrPic}(\Rep(\ku G))$ for a finite Abelian group $G$. This group is isomorphic
to the group of automorphism of $G\oplus \widehat{G}$, here $\widehat{G}$ is
the group of characters of $G$, such that they preserve the quadratic form $q:G\oplus \widehat{G}
\to \ku$, $q(g,f)=f(g)$. In this section we use the same ideas to give a more
compact description of the group $\nic(V,u,G)$.

\medbreak

Let
$(W,\beta, \alpha)\in \ere(V, u,G)$. 
Set $\tau(W,\beta)$  the subspace of $V\oplus V^*\oplus V\oplus V^*$ defined by
$$\{(w_1,f_1,w_2,f_2): (w_1,w_2)\in W, \,
(f_1,f_2)\in W^*,\,
\widehat{\beta}(w_1,w_2)=f_1-f_2\}.  $$
Recall the definition of $\widehat{\beta}$ given in \eqref{bilinear-f1}. If
$(W',\beta', \alpha')$ is another element in $\ere(V, u,G)$ we denote
$\tau(W,\beta)\bullet \tau(W',\beta')$ the set of elements 
$(w_1,f_1,w_2,f_2)$ such that there exists a unique $(v,g)\in V\oplus V^*$ such that
$(w_1,f_1,v,g)\in \tau(W,\beta)$ and $(v,g,w_2,f_2)\in \tau(W',\beta')$.

\medbreak
Let
$\lag(V, u,G)$ be the set of pairs $( \tau(W,\beta), \alpha)$ where
$(W,\beta, \alpha)$ is an invertible element in $\ere(V, u,G)$. If $( \tau(W,\beta),\alpha), ( \tau(W',\beta'), \alpha')\in \lag(V, u,G)$
define 
\begin{equation}\label{g-in-lag}( \tau(W,\beta), \alpha)\bullet ( \tau(W',\beta'), \alpha')= (\tau(W,\beta)\bullet  \tau(W',\beta'),
\alpha\circ\alpha').
\end{equation}
Two elements $( \tau(W,\beta), \alpha), ( \tau(W',\beta'), \alpha')$ in $\lag(V, u,G)$ are \emph{equivalent}  
if there exists $(x,y)\in G\times G$ such that 
$$( \tau(W',\beta'), \alpha')=( \tau((x,y)\cdot W,(x,y)\cdot \beta), \alpha).$$ We denote by $\lagb(V, u,G)$
 the set of equivalence classes in $\lag(V, u,G)$. The next lemma is an analogue result
of \cite[Prop. 10.3]{ENO}. 
\begin{lema} The set  $\lagb(V, u,G)$ is a group with operation defined by 
\eqref{g-in-lag} in each equivalence class and identity element the class of
$(\{(v,f,v,f): v\in V, f\in V^*\}, \id)$. The map $\tau: \nic(V, u,G)\to
\lagb(V, u,G)$ that sends the class of $(W,\beta, \alpha)$ to
the class of $( \tau(W,\beta), \alpha)$ is a group
isomorphism.
\end{lema}
\pf  The proof that  $\lagb(V, u,G)$ is a group is straightforward. Let us take $(W,\beta,\alpha)$, $(W',\beta', \alpha')\in  \ere(V, u,G)$
and $(w_1,f_1,w_2,f_2)\in \tau(W,\beta)\bullet  \tau(W',\beta')$. Then
there exists $(v,g)\in V\oplus V^*$ such that
$(w_1,f_1,v,g)\in \tau(W,\beta)$ and $(v,g,w_2,f_2)\in \tau(W',\beta')$. Hence
$$\widehat{\beta}(w_1,v)=f_1-g, \quad \widehat{\beta'}(v,w_2)=g-f_2,$$
which implies that 
$$\widehat{\beta\bullet \beta'}(w_1,w_2)=f_1-f_2.$$
Thus, $(w_1,f_1,w_2,f_2)\in \tau(W \bullet  W', \beta\bullet \beta')$ and we have 
an inclusion $ \tau(W,\beta)\bullet  \tau(W',\beta')\subseteq \tau(W \bullet  W', \beta\bullet \beta')$.
The other inclusion is proven similarly. Thus $\tau$ is  well-defined and injective.
By definition of $\lag(V, u,G)$ the map $\tau$ is surjective.
\epf

 The group $G\times G$ acts on the set of linear maps $T:V\oplus V^*\to V\oplus V^*$
as follows. If $(x,y)\in G\times G$, $(v,f)\in V\oplus V^*$ define
\begin{equation}\label{defi-lag}( (x,y)\cdot T)(v,f)=x^{-1}\cdot T(y\cdot v,y\cdot f). \quad
\end{equation}
The action of $G$ on $V^*$ is given by
$$(x\cdot f)(v)= f(x^{-1}\cdot v),$$
for all $x\in G$, $f\in V^*$, $v\in V$.
\begin{defi}\label{o-group} Let $\oc(V,u,G)$  the set of pairs $(T,\alpha)$ 
where 
\begin{itemize}
 \item[(i)] $\alpha\in O(G\oplus \widehat{G})$ such that $(u,u)\in U_\alpha$,
 \item[(ii)] $T:V\oplus V^*\to V\oplus V^*$ is a linear isomorphism such that
\begin{equation}\label{eqt-1} (x,y)\cdot T=T, \quad \text{ for all }(x,y)\in U_\alpha,
\end{equation}
\begin{equation}\label{eqt-2} T^1(0,f)=0, \quad T^2(0,f)(T^1(v,0))=f(v), \quad \text{ for all }
f\in V^*, v\in V.
\end{equation}
Here $T(v,f)=(T^1(v,f), T^2(v,f))$ for all $f\in V^*, v\in V$.
\end{itemize}

\end{defi}

Two elements $(T,\alpha)$ , $(T',\alpha')$ are \emph{equivalent} if there exists
 $(x,y)\in G\times G$ such that 
$$T'=(x^{-1},y^{-1})\cdot T, \;\; \alpha=\alpha'.$$
The class of an element $(T,\alpha)\in \oc(V,u,G)$ will be denoted by
$\overline{(T,\alpha)}$ and the set of equivalence classes will be denoted $\ocb(V,u,G)$.

\begin{rmk} If $(T,\id)\in \oc(V,u,G)$ then $T\in \Aut_G(V\oplus V^*)$.
\end{rmk}

\begin{lema} The set $\ocb(V,u,G)$ is a group with unit element 
$\overline{(\Id,\id) }$ and composition 
$$\overline{(T,\alpha)}\bullet \overline{(T',\alpha')}=
 \overline{(T\circ T',\alpha\circ \alpha')},$$
for all $\overline{(T,\alpha)}, \overline{(T',\alpha')}\in \ocb(V,u,G)$.\qed
\end{lema}

\begin{teo} There is an isomorphism of groups $\nic(V,u,G)\simeq\ocb(V,u,G)$.
\end{teo}
\pf Let $(T,\alpha)$ be a representative of a class in $\ocb(V,u,G)$.
Define $T^1:V\oplus V^*\to V, T^2:V\oplus V^*\to V^*$ by
$T(v,f)=(T^1(v,f), T^2(v,f))$ for any $(v,f)\in V\oplus V^*$. Let
$W_T$ the subspace of $V\oplus V$ defined as
$$W_T=\{(T^1(v,f),v): v\in V, f\in V^*\},$$
and the bilinear form $\beta_T:W_T\times W_T\to \ku$ defined by
$$\beta_T((T^1(v_1,f_1),v_1), (T^1(v_2,f_2),v_2))=T^2(v_1,f_1)(T^1(v_2,f_2))- f_1(v_2),$$
for all $(v_1,f_1), (v_2,f_2)\in  V\oplus V^*$. 
\begin{claim} $(W_T, \beta_T, \alpha)\in  \ere(V, u,G)$.
\end{claim}
\pf[Proof of Claim] Let us prove that $ \beta_T$ is $U_\alpha$-invariant.
The other conditions can be easily verified. Let $(g,h)\in U_\alpha$, $(v_1,f_1), (v_2,f_2)
\in V\oplus V^*$ then
$\beta_T((g,h)\cdot (T^1(v_1,f_1),v_1), (g,h)\cdot (T^1(v_2,f_2),v_2))$ is equal to
\begin{align*}&=\beta_T( (g\cdot T^1(v_1,f_1),h \cdot v_1),  (g\cdot T^1(v_2,f_2),h \cdot v_2))\\
&=\beta_T( (T^1(h \cdot v_1,h \cdot f_1),h \cdot v_1),  ( T^1(h \cdot v_2,h \cdot f_2),h \cdot v_2))\\
&= T^2(h \cdot v_1,h \cdot f_1)(T^1(h \cdot v_2,h \cdot f_2))- h \cdot f_1(h \cdot v_2)\\
&=T^2(v_1,f_1)(T^1(v_2,f_2))- f_1(v_2)\\
&=\beta_T((T^1(v_1,f_1),v_1), (T^1(v_2,f_2),v_2)).
\end{align*}
The second and fourth equalities follows because $(g,h)\cdot T=T$.\epf

We will establish an isomorphism $\sigma:\ocb(V,u,G)\to\lagb(V, u,G)$ defined by
$$\sigma\overline{(T,\alpha)}=\overline{(\tau(W_T, \beta_T),\alpha)},$$
for all $\overline{(T,\alpha)}\in \ocb(V,u,G).$ This map does not depend on the representative
class of $\overline{(T,\alpha)}$. Let us prove that it is injective. Let
$\overline{(T,\alpha)}\in \ocb(V,u,G)$ such that $$\overline{(\tau(W_T, \beta_T), \alpha)}
=\overline{ (\{(v,f,v,f): v\in V, f\in V^*\}, \id)}.$$
Since $ (\{(v,f,v,f): v\in V, f\in V^*\}=\tau(\diag(V),0)$ there exists an element
$(x,y)\in G\times G$ such that
\begin{align*}U_\alpha= \diag(G) ,\;\; \psi_\alpha=1, \;\; W_T=\{(x\cdot v, y\cdot v): v\in V\}, \;\;
\beta_T=0.
\end{align*}
This implies that $T^1(v,f)=xy^{-1}\cdot v$ for all $(v,f)\in V\oplus V^*$ and since
$\beta_T=0$ then $T^2(v,f)=xy^{-1}\cdot f$, thus $T=(x,y)^{-1}\cdot \Id$. 
Hence $\overline{(T,\alpha)}=\overline{(\Id,\id)}$ and $\sigma$ is
injective. Finally, let us prove that $\sigma$ is surjective. Let 
$ (\tau(W,\beta),\alpha)\in \lag(V, u,G)$. If $(w_1,f_1,w_2,f_2), 
(w'_1,f'_1,w_2,f_2)\in \tau(W,\beta)$ then $(w_1-w'_1,0)\in W$
which implies that $w_1=w'_1$. Also
$$\widehat{\beta}(w_1,w_2)=f_1-f_2=f'_1-f_2,$$
thus $f'_1=f_1$. In conclusion the pair $(w_1,f_1)$ depends on 
$(w_2,f_2)$, therefore there is a linear function 
$T:V\oplus V^*\to V\oplus V^*$ such that
$W=W_T$. If the element $(0,0,v,f)\in W$ then $v=0$, $f=0$, thus
$T$ must be injective and consequently bijective. It is not difficult to
see that $\beta=\beta_T$. This finishes the proof that $\sigma$ is surjective
and the proof of the Theorem.
\epf

\begin{exa} Suppouse $\ku=\C$. Let $\Z_2$ be the cyclic group of order 2 with generator $u$. Let
$V$ be a finite-dimensional vector space such that $u$ acts as $-1$ on $V$. Set 
$H=\wedge(V)\#\ku \Z_2$. Assume $\dim V=1$, so $H$ is the Sweedler's Hopf algebra.

The group $O(\Z_2\oplus \widehat{\Z_2})=\{\id, \gamma\}$, see example \ref{exa-brC2}.
Note that $ U_\gamma=\Z_2\oplus  \Z_2$. Define
$$ \oc=\{ A\in SL_2(\C): A_{12}=0\}. $$
The Brauer-Picard group of $\Rep(H)$ is isomorphic to the group $\oc\times \Z_2$. In particular
for any $\xi\in \ku $ the matrices 
$$\left(\begin{matrix} i \\  \xi
\end{matrix}
\begin{matrix}\;\;  0\\  \; \;\; -i \end{matrix}\right)$$
give a one parameter family  invertible bimodule categories over $\Rep(H)$ of order 4. 
\end{exa}

\bibliographystyle{amsalpha}

\end{document}